\documentclass[10pt]{article}
\usepackage{amssymb,amsmath,latexsym}
\input xypic

\newtheorem{theorem}{Theorem}[subsection]
\newtheorem{proposition}[theorem]{Proposition}
\newtheorem{corollary}[theorem]{Corollary}
\newtheorem{lemma}[theorem]{Lemma}

\newcommand{\qed}{ \hfill $\blacksquare$ \bigskip}
\newcommand{\eb}{\qed}
\newcommand{\bew}{{\em Proof }}

\newcommand{\bG}{{\bf G}}
\newcommand{\RR}{\mathbb{R}}
\newcommand{\Q}{\mathbb{Q}}
\newcommand{\Z}{\mathbb{Z}}
\newcommand{\Nat}{\mathbb{N}}
\newcommand{\A}{\mathcal{A}}
\newcommand{\emp}{\,\mathcal{E}\,}
\renewcommand{\P}{\mathcal{P}}

\newcommand{\C}{\mathcal{C}}
\newcommand{\D}{\mathcal{D}}
\newcommand{\T}{\mathcal{T}}

\newcommand{\erz}[1]{<\!#1\!>}
\newcommand{\im}{\hbox{\rm im}}
\newcommand{\s}{{\bf d}}
\newcommand{\f}{{\bf r}}
\newcommand{\E}{\, \mathcal{E} \,}
\newcommand{\La}{\D}

\title{Universal groups for point-sets and tilings}
\author{Johannes Kellendonk\\
School of Mathematics\\
University of Wales, Cardiff\\
Cardiff CF24 4YH\\
United Kingdom
\and
Mark~V.~Lawson\\ 
Mathematics Divison\\ 
 School of Informatics\\
University of Wales, Bangor\\
Dean Street\\
Bangor, Gwynedd LL57 1UT\\
United Kingdom}

\begin{document}
\maketitle
\begin{abstract} 
\noindent
We study the universal groups of inverse semigroups associated with 
point sets and with tilings.
We focus our attention on two classes of examples.
The first class consists of point sets which are obtained
by a cut and projection scheme (so-called model sets). Here we introduce
another inverse semigroup which is given in terms of the defining
data of the projection scheme and related to the model set by the empire
congruence.
The second class is given by one-dimensional tilings.
\end{abstract}
\medskip

\section{Introduction}
This paper concerns the algebraic description of aperiodic
systems which arise in solid state physics. 
These systems are point-sets, typically Delone sets, or tilings in $\RR^d$. 
To describe such systems we use inverse semigroups but
these should be viewed as stepping-stones to our real goal which,
as the title suggests, is to associate groups with them.

The question of what we mean by an ordered,
as opposed to an unordered, solid 
is a difficult one and attempts 
to solve it draw on many different
areas of mathematics and physics; 
see \cite{BaakeMoody00} for an overview of the latest developments. 
The usual mathematical starting point is the point-set of 
equilibrium positions of atoms in a solid, or of a tiling which
encodes these atomic positions. 
The traditional approach to order in solids is based on their ($X$-ray)
diffraction images and,
until about twenty years ago, 
this was modelled 
by periodic structures.\footnote{A point-set describing such a structure  
is given by (finitely many) orbits of a lattice and
hence is largely characterized by a group, 
its crystallographic group,
which is a semidirect product of the lattice with the point group.} 
But in recent years, particularly with the discovery of
quasicrystals, it has become clear that ordered structures can be
aperiodic.

In this context, Lagarias proposed a hierarchy of point-sets
\cite{Lagarias} to describe various levels of order. The highest
level of aperiodic order furnish the so-called model sets:
these are sets described by a cut-and-projection scheme,
and are used to describe 
quasicrystals; we discuss them in some detail in this paper.
In Lagarias' characterization of a Delone set $\D$ the properties of
the set of difference vectors $\D-\D$ play an important role.
In particular, 
he considers the subgroup $H_\D$ of $\RR^d$ generated by $\D-\D$. This
group is important for questions related to (X-ray) diffraction of the
material 
(but we will have nothing to say about this here). 
In this paper, we will
construct an {\em a priori} different group $G_\D$ from $\D$ 
and compare it with $H_\D$. 
As we shall show,
$G_\D$ is the universal group of the point-set semigroup of $\D$; 
this point-set semigroup can be defined for any point-set of $\RR^d$
by a construction which resembles the construction of the
tiling semigroup from a tiling \cite{KL2000}.
In general, however, $G_\D$ is not abelian.
The property of its being finitely generated can serve 
as the basis
for a characterization of a potentially new class of Delone sets in 
Lagarias' hierarchy. 

We discuss the inverse semigroup of a model set and its universal
group at some length. In fact, for model sets we can 
define another inverse semigroup which is derived  from
the data which enter into the cut-and-projection scheme defining the model set. 
This allows us to prove that
for a large class of model sets $G_\D$ is isomorphic to $H_\D$. 
On the other hand, we point out that not all model sets
have universal groups isomorphic to $H_\D$.

The first structural descriptions of quasicrystals used tilings not 
point-sets.
This has advantages if one wants to describe self-similar structures
i.e.\ structures which possess an inflation/deflation symmetry, 
a process which involves sub-dividing tiles. 
The Penrose tilings are the most prominent examples of this. 
Also a recent idea explaining the stability of pentagonal
quasicrystals uses an extension of the tiling idea to
overlapping clusters \cite{SteinhardtJeong96} (Gummelt's decagons
\cite{Gummelt96}) 
something which mimics the multiplication law in tiling semigroups. 
On the other hand, for
most of the mathematical theory on aperiodic structures, in particular
the topological dynamical aspects and parts of diffraction theory, it does not
matter whether one works with a tiling or a point-set description
because they depend
only on a topological groupoid derived from the
semigroup and this groupoid is  (up to equivalence)
the same whether one starts with a
tiling or a  point-set semigroup.
The universal group of a tiling semigroup, however,
differs quite substantially from that of a Delone set derived from the
tiling. We will show below that the universal group 
for a one-dimensional tiling is always the free group on as many
generators as there are (equivalence classes) of touching pairs of tiles.
This difference from the point-set case 
comes from the fact that the inverse semigroups associated to a
tiling or a Delone set contain different information about the local
structure. 
We have not yet fully understood the significance of this.

The main result of the paper is Theorem~3.4.5 of Section~3
which describes, given certain assumptions,
the universal groups associated with point-sets
obtained by the cut-and-projection scheme.
A key ingredient in the proof of this result is
a theorem due to Macbeath concerning group presentations
\cite{Macbeath64}.
Section~4 deals briefly with the universal groups of
connected tiling semigroups for one-dimensional tilings
--- it is intended as a contrast to the results 
obtained in Section~3.
Section~2 provides all the algebraic preliminaries.\\

\noindent
{\bf Acknowledgements} The second author is grateful to Stuart Margolis
for conversations concerning the universal group of a semigroup,
and for pointing out reference \cite{Nica94} which suggested the approach
adopted in Section~2.

\section{Background}\setcounter{theorem}{0}
The aim of this paper is to find ways of associating groups
with certain mathematical models of quasicrystals.
Our method is to associate inverse semigroups with the different models
of such quasicrystals and then construct the universal group of the inverse
semigroup.
In this section, we outline the algebraic results needed to
read this paper.

\subsection{Inverse semigroups}
We shall principally be concerned with sets $S$ equipped
with partial binary operations $\circ$.
We write $\exists s \circ t$ to mean that $s \circ t$ is defined.
Usually we shall denote the partial binary operation
by concatenation.

Let $(S,\circ)$ and $(T,\ast)$ be two sets equipped with 
partial binary operations.
A function $\theta \colon \: S \rightarrow T$ 
is called a {\em morphism}
if $\exists s \circ t$ in $S$ implies that
$\exists \theta (s) \ast \theta (t)$ and
$\theta (s \circ t) = \theta (s ) \ast \theta (t)$.  
A morphism is called
a {\em homomorphism} if it satisfies the additional property
that $s \circ t$ exists if and only if
$\theta (s) \ast \theta (t)$ exists.
An {\em isomorphism} is a bijective homomorphism.\footnote{It is
easy to check that its inverse is also a homomorphism.}

A {\em group-like set} \cite{Jekel77} is a pair $(S,\circ)$ consisting
of a set $S$ equipped with a partial binary operation $\circ$ 
satisfying the following three axioms:
\begin{description}

\item[{\rm (GL1)}] There is an element $1 \in S$ such
that $\exists 1 \circ s$ and $\exists s \circ 1$ 
for all $s \in S$ and
$1 \circ s = s = s \circ 1$.

\item[{\rm (GL2)}] For each $s \in S$ there exists an element
$s^{-1} \in S$ such that $\exists s^{-1} \circ s$ and $\exists s \circ s^{-1}$
and both products are equal to $1$.

\item[{\rm (GL3)}] If $\exists s \circ t = u$ then
$\exists t^{-1} \circ s^{-1}$ and equals $u^{-1}$.
\end{description}
Groups are clearly special cases of group-like sets
as are the pregroups of Stallings \cite{Stallings71};
however, Jekel gives an interesting geometrically motivated group-like set
which is not a pregroup in \cite{Jekel76}.
Two group-like sets $(S,\circ)$ and $(T,\ast)$ are said to
be {\em isomorphic} if there is an isomorphism
in the sense above
$\alpha \colon \: S \rightarrow T$
which in addition satisfies 
$\alpha (1) = 1$,
and 
$\alpha (s^{-1}) = \alpha (s)^{-1}$.

Often the partial binary operation will have the additional
property that $(s \circ t) \circ u$
exists if and only if $s \circ (t \circ u)$ exists
in which case they are equal.
We call such a structure a {\em semigroup}.\\

\noindent
{\bf Remark} A semigroup
is usually required to have an everywhere defined operation.
However, a semigroup in our sense becomes a semigroup
in the usual sense by simply adjoining a zero
and defining all undefined products to be zero.\\

A {\em monoid} is a semigroup $S$ with a distinguished element
$1$ such that $\exists 1s$ and $\exists s1$ for all $s \in S$
and $1s = s = s1$.
An important example of a monoid is the following.
Let $\Sigma$ be an alphabet; that is, a non-empty set.
The {\em free monoid} on $\Sigma$, denoted by $\Sigma^{\ast}$,
consists of all strings over $\Sigma$ with
concatenation as product and identity the empty string $\epsilon$.
An {\em idempotent} in a semigroup
is an element $e \in S$ such that
$\exists e \circ e$ and $e \circ e = e$.
A semigroup is itself {\em idempotent} if every element is idempotent.
A semigroup is {\em commutative} if
$\exists s \circ t$ if and only if $\exists t \circ s$
in which case they are equal.  
We refer the reader to \cite{Howie} for standard results
about semigroups.

A semigroup $S$ is said to be {\em inverse}
if for each $s \in S$ there exists a unique element
$s^{-1} \in S$ such that 
$ss^{-1}s$ and $s^{-1}ss^{-1}$ are both defined and
$s = ss^{-1}s$ and $s^{-1} = s^{-1}ss^{-1}$.
The theory of inverse semigroups is developed in \cite{Lawson98}.
We record two important results here that we will need later.
The set of idempotents $E(S)$ of an inverse semigroup $S$ forms
a commutative idempotent semigroup or {\em semilattice}.
On each inverse semigroup $S$, 
we may define a partial order, called the
{\em natural partial order}, 
by $s \leq t$ if and only if there is an idempotent $e$ such that $s = te$.
Morphisms of inverse semigroups preserve the natural partial order:
if $\theta \colon \: S \rightarrow T$ is a morphism
of inverse semigroups and $s \leq t$ then
$s = te$ and so $\theta (s) = \theta (te) = \theta (t)\theta(e)$
but $\theta (e^{2}) = \theta (e)$.
Thus $\theta (s) \leq \theta (t)$.
If $\theta \colon S \rightarrow G$ is a morphism 
from an inverse semigroup to a group 
then for each idempotent $e \in S$ we have that
$\theta (e) = 1$, and so
$s \leq t$ implies that $\theta (s) = \theta (t)$.
A morphism $\theta \colon \: S \rightarrow T$ 
between inverse semigroups
is said to be {\em idempotent pure}
if $\theta (a)$ an idempotent implies that
$a$ is an idempotent.

An inverse semigroup $S$ is said to
be {\em strongly $E^{\ast}$-unitary}
if there is an idempotent pure morphism from $S$ to a group.
An inverse semigroup is {\em $F^{\ast}$-inverse}
if every element is beneath a unique maximal element
with respect to the natural partial order. 
An $F^{\ast}$-inverse semigroup which is also
strongly $E^{\ast}$-unitary
is said to be {\em strongly $F^{\ast}$-inverse}.

A {\em category} $C$ is a semigroup (in the terms of this paper)
satisfying certain extra conditions:
there is a subset $C_{o}$ of $C$ called the set
of {\em identities} which are idempotents $e$ such that
whenever $xe$ or $ex$ are defined they equal $x$;
there are maps $\s,\f \colon \C \rightarrow C_{o}$ such that
$\exists x \s (x)$ and $\exists \f (x) x$ 
and $\exists xy$ if and only if $\s (x) = \f (y)$;
in which case $\s (xy) = \s (y)$ and $\f (xy) = \f (x)$.
It is best to think of a category as a directed
graph in which each element $x$ is an {\em arrow}
from $\s (x)$ to $\f (x)$.
We call $\s (x)$ the {\em right identity} of $x$
and $\f (x)$ the {\em left identity} of $x$.
If $x,y \in C$ are such that
$\s (x) = \s (y)$ and $\f (x) = \f (y)$
then we say $x$ and $y$ are {\em parallel}.
An element $x \in C$ such that
$\s (x) = \f (x)$ is called a {\em loop}.
If $e \in C_{o}$ then the set of
all $x \in C$ such that 
$\s (x) = \f (x)$ forms a monoid 
called the {\em local monoid (at $e$)}.
A category is {\em locally idempotent (resp. locally commutative)}
if all the local monoids are idempotent (resp. commutative).
A category whose semigroup is inverse 
is called {\em inverse}.
In an inverse category,
if $s \leq t$ then $s$ and $t$ must be parallel.
The most important morphisms between categories
are the functors:
$\theta \colon \: C \rightarrow D$ is a {\em functor}
if $\theta$ maps identities to identities,
$\theta (\s (x)) = \s (\theta (x))$,
$\theta (\f (x)) = \f (\theta (x))$,
and if $xy$ is defined then $\theta (x) \theta (y)$ is
defined and $\theta (xy) = \theta (x) \theta (y)$.

We now describe a procedure for
constructing examples of strongly $E^{\ast}$-unitary inverse semigroups;
all the inverse semigroups occurring in this paper will be of this type. 
More details can be found in \cite{KL2000}.

Let $G$ be a group and $X$ a set.
We say that $G$ {\em acts partially on $X$ }
if there is a partial function $G \times X \rightarrow X$
denoted by $(a, x) \mapsto a \cdot x$ satisfying the
following three axioms:

\begin{description}
\item[{\rm (PA1)}] $\exists 1 \cdot x$ for each $x \in X$ and $1 \cdot x = x$.

\item[{\rm (PA2)}] $\exists g \cdot (h \cdot x)$ implies $\exists (gh) \cdot x$
and $g \cdot (h \cdot x) = (gh) \cdot x$.

\item[{\rm (PA3)}] $\exists g \cdot x$ implies that 
$\exists g^{-1} \cdot (g \cdot x)$, and $g^{-1} \cdot (g \cdot x) = x$.
\end{description}

We say that the action is {\em free}
if $\exists g \cdot x = x$ implies that $g = 1$.

The following lemma will be useful below.

\begin{lemma} Let the group $G$ act partially on the set $X$.
Suppose that $\exists b \cdot x$ and $\exists (ab)\cdot x$
then $\exists a \cdot (b \cdot x)$. 
\end{lemma}
\bew\ By assumption $\exists (ab) \cdot x$.
Thus by (PA3), $\exists (ab)^{-1} \cdot ((ab) \cdot x)$ and is equal to $x$.
By assumption $\exists b \cdot x$.
Thus $\exists b \cdot ((ab)^{-1} \cdot ((ab) \cdot x))$.
By (PA2), we have that
$\exists (b(ab)^{-1}) \cdot ((ab) \cdot x))$.
Thus
$\exists a^{-1} \cdot ((ab) \cdot x)$.
By (PA3) we have that
$\exists a \cdot (a^{-1} \cdot ((ab) \cdot x))$.
Thus by (PA2)
$\exists a \cdot (b \cdot x)$.\qed

Let $S$ be a category whose partial product we denote by
concatenation and let $G$ be a group.
We say that $G$ {\em acts freely and partially on} $S$
if $G$ acts partially on both $S_{o}$
and $S$, and freely on $S_{o}$;
in addition, the following four axioms are satisfied:
\begin{description} 
\item[{\rm (A1)}] If $\exists g \cdot x$ then
$\exists g \cdot \s (x)$ and $\exists g \cdot \f (x)$ and
$\s (g \cdot x) = g \cdot \s (x)$ and
$\f (g \cdot x) = g \cdot \f (x)$.

\item[{\rm (A2)}] If $\exists xy$ and $\exists g \cdot x$ and
$\exists g \cdot y$ then $\exists g \cdot (xy)$
and $g \cdot (xy) = (g \cdot x)(g \cdot y)$.
(Observe that $\exists (g \cdot x) (g \cdot y)$ by (A1)).

 \item[{\rm (A3)}] If $\exists g \cdot (xy)$ then $\exists g \cdot x$
and $\exists g \cdot y$.

\item[{\rm (A4)}] If $\exists g \cdot (xy)$ 
then $\exists (g \cdot x)(g \cdot y)$
and $g \cdot (xy) = (g \cdot x)(g \cdot y)$.
(Observe that $\exists g \cdot x$ and $\exists g \cdot y$ by (A3)).
\end{description}

The assertion in part~(ii) of the following result
that the semigroup is `strongly $E^{\ast}$-unitary'
is due to John Fountain
(private communication).\footnote{By replacing the categories
in (ii) with `multiplicative graphs' we can actually construct
every strongly $E^{\ast}$-unitary inverse semigroup in this way;
this is a reformulation of the results in \cite{Steinberg2000}.}

\begin{theorem} Let $G$ be a group which acts freely and partially
on the category $S$.
\begin{description}
\item[{\rm (i)}] Define $\sim$ on $S$ by $x \sim y$ iff there exists $g \in G$
such that $g \cdot x = y$.
Then $\sim$ is an equivalence relation.
Let $S/G = \{[x] \colon \: x \in S \}$ be the set of equivalence classes.
Define the following operation on $S/G$:
$$[x][y] = [(g \cdot x)(h \cdot y)]$$
if there exist $g,h \in G$ such that
$(g \cdot x)(h \cdot y)$ is defined in $S$.
Then with respect to this operation,
$S/G$ is a semigroup.
If the category is inverse then so too is $S/G$.

\item[{\rm (ii)}] Suppose $S$ is an inverse category
whose local monoids are idempotent.
Then $S/G$ is a strongly $E^{\ast}$-unitary inverse semigroup.
If $S$ has the further property that each set of parallel elements has a 
maximum then $S/G$ is strongly $F^{\ast}$-inverse.\qed
\end{description}
\end{theorem}

Theorem~2.1.2 can be used to construct
inverse semigroups from tilings and point-sets.
Here are some examples.\\

\noindent
{\bf Examples}
\begin{description}

\item[{\rm (i)}] {\em Tiling semigroups} 

Let $\T$ be a tiling in $\mathbb{R}^{n}$ which we view as a countable
set of bounded closed subsets, called the tiles; see \cite{KL2000} for
the detailed requirements. 
The set of triples $(a,P,b)$ where $P$ is any finite subset
of $\T$ and $a$ and $b$ belong to $P$ forms
an inverse category $C$ when we define a partial product
by $(a,P,b)(b,Q,c) = (a, P \cup Q,c)$ and undefined
otherwise. The identities are the elements of
the form $(a,a,a)$, and $(a,P,b)^{-1} = (b,P,a)$.
We have that $(a,P,b) \leq (c,Q,d)$ if and only if
$a = c$ and $b = d$ and $Q \subseteq P$.
The local monoids of $C$ consist entirely
of idempotents and each set of parallel elements
contains a maximum element.
Let $G$ be the group of translations of $\mathcal{R}^{n}$.
This acts on $C$ freely and partially.
We denote the resulting strongly $F^{\ast}$-inverse semigroup by $\Gamma(\T)$.
It is the original {\em tiling semigroup} considered 
by one of the authors \cite{Kellendonk97a}.

The tiling semigroup $\Gamma(\T)$ contains a subsemigroup $S(\T)$
which is given by equivalence classes of triples $(a,P,b)$ where $P$
is connected in the sense that the union of its tiles cover a
connected subset of $\RR^d$. We call $S(\T)$ here 
the {\em connected tiling semigroup} of $\T$.
It is strongly $E^{\ast}$-unitary but not $F^{\ast}$-inverse
in general, except in the case of 1-dimensional tilings,
which we examine in Section~4. 

\item[{\rm (ii)}] {\em Point-set semigroups}

Point-set semigroups are defined in much the same way as tiling
semigroups. 
Let $\La \subseteq \mathbb{R}^{n}$.
The set of triples $(a,P,b)$ where $P$ is a finite subset
of $\La$ and $a,b \in P$ forms an inverse category where
the product and inverses are defined analogous to that in (i).
Once again we have a free partial action of the group of translations
and the quotient of the category with respect to this action yields
a strongly $F^{\ast}$-inverse semigroup.
We denote this inverse semigroup by $\Gamma (\La)$
and call it the {\em point-set semigroup} of $\La$.

\item[{\rm (iii)}] {\em The semigroups $\Gamma (X,G,H)$} 

Here is a different class of examples which we mainly introduce to
describe point-set semigroups of model sets more effectively.
Let $H$ be a group,
let $G$ be a subgroup of $H$
and let $X$ be a subset of $H$ containing the identity.
Let $\A = \{gX \colon \: g \in G \mbox{ and } 1 \in gX\}$,
and let $\P$ be the set of all finite intersections
of elements of $\A$;
in other words, $\P$ is the semilattice generated by
$\A$ under intersection.
Let $C$ be the set of all triples
$(a,P,b)$ where $P \in \P$ and $P \subseteq aX, bX$.
Define $(a,P,b)(b,Q,c) = (a, P \cap Q,c)$,
which is easily seen to be well-defined.
Observe that $(b,bX,b)$ is a right identity
of $(a,P,b)$ and $(a,aX,a)$ is a left identity.
We define $(a,P,b)^{-1} = (b,P,a)$.
It is now easy to check that $C$ is an inverse category.
The group $G$ acts partially on $C$
by $g \cdot (a,P,b) = (ga, gP, gb)$ if $gP \in \P$.
It is easy to check that $G$ acts freely.
The inverse semigroup $C/G$ is denoted by $\Gamma (X,G,H)$.
\end{description}

For particular choices of $X,G,H$ the last semigroup is strongly related to
point-set semigroups for model sets. This relation is based on the
notion of the `empire' of a tile or pattern in a tiling which can be found  
in the book of Gr\"unbaum and Shephard \cite{GS}.
Zhu \cite{Zhu} has shown that this notion defines a congruence on tiling
semigroups.
Zhu's definition can easily be extended to any inverse semigroup
constructed from a group acting freely and partially
on a locally idempotent inverse category.
We now show how this can be done.

Let $G$ be a group acting freely and partially on
the locally idempotent inverse category $S$.
Define the relation $\equiv$ on $S$
by $x \equiv y$ iff $x$ and $y$ are parallel and 
for each $g \in G$ we have that
$\exists g \cdot x \Leftrightarrow \exists g \cdot y$.
Now define the relation $\E$ on $S/G$ as follows:
$[x] \E [y]$ iff there exist $x' \in [x]$ and $y' \in [y]$
such that $x \equiv y$.

\begin{proposition} Let $G$ be a group acting freely and
partially on the locally idempotent inverse category $S$.
Then the relation $\E$ is an idempotent pure congruence on $S/G$.
\end{proposition}
\bew\ We show first that $\E$ is an equivalence relation.
It is clear that $\E$ is reflexive and symmetric,
so we need only prove that it is transitive.

Suppose that $[x] \E [y]$ and $[y] \E [z]$.
Let $a \cdot x \in [x]$ and $b \cdot y \in [y]$ be such that
$a \cdot x \E b \cdot y$ and
let $c \cdot y \in [y]$ and $d \cdot z \in [z]$ be such that
$c \cdot y \equiv d \cdot z$.
First $\exists b \cdot ( c^{-1} \cdot (c \cdot y))$ by (PA3).
Thus $(bc^{-1}) \cdot (c \cdot y)$ by (PA2).
But $c \cdot y \equiv d \cdot z$.
Thus $\exists (bc^{-1}) \cdot (d \cdot z)$.
Hence $\exists (bc^{-1}d) \cdot z$.
By (A1), we have that
$$\s ((bc^{-1}d) \cdot z) 
= 
\s ((bc^{-1}) \cdot (d \cdot z))
=
(bc^{-1}) \cdot \s (d \cdot z)
=
(bc^{-1}) \cdot \s (c \cdot y)
=
(bc^{-1}) \cdot (c \cdot \s (y))
=
b \cdot \s (y)
=
\s (b \cdot y)
=
\s (a \cdot x).$$
Similarly
$\f ((bc^{-1}d) \cdot z) 
=
\f (a \cdot x)$.
Thus $a \cdot x$ and $(bc^{-1}d) \cdot z$ are parallel.
Suppose that $\exists g \cdot (a \cdot x)$.
Then $\exists g \cdot (b \cdot y)$ since
$a \cdot x \equiv b \cdot y$.
But $y = c^{-1} \cdot (c \cdot y)$.
Thus
$\exists (gbc^{-1}) \cdot (cy)$.
It follows that
$\exists (gbc^{-1}) \cdot (d \cdot z)$ since
$c \cdot y \equiv d \cdot z$.
Thus $\exists (gbc^{-1}d) \cdot z$ and $\exists (bc^{-1}d) \cdot z$
so $\exists g \cdot ((bc^{-1}d) \cdot z)$ by Lemma~2.1.1.
Conversely, suppose that
$\exists g \cdot ((bc^{-1}d) \cdot z)$.
Then $\exists (gbc^{-1}d) \cdot z$ and $d \cdot z$
and so by Lemma~2.1.1, we have that
$(gbc^{-1}) \cdot (d \cdot z)$.
Thus $(gbc^{-1}) \cdot (c \cdot y)$ since
$c \cdot y \equiv d \cdot z$.
Thus $\exists (gb) \cdot y$ and $\exists b \cdot y$
so $\exists g \cdot (b \cdot y)$ by Lemma~2.1.1.
Thus $g \cdot (a \cdot x)$ since
$b \cdot y \equiv a \cdot x$.
We have therefore proved that $a \cdot x \equiv (bc^{-1}d) \cdot z$
and so $[x] \E [z]$, as required.

We now have to show that the equivalence relation $\E$
is actually a congruence on the semigroup $S/G$.
Let $[x] \E [y]$ and $[u] \E [v]$ and suppose that
$[x][u]$ and $[y][v]$ are both defined.
We prove that $[x][u] \E [y][v]$. 
Suppose that $[x][u] = [(a \cdot x)(b \cdot u)]$,
$m \cdot x \equiv n \cdot y$ and $p \cdot u \equiv q \cdot v$.
It is easy to check that
$a \cdot x$ and $(am^{-1}n) \cdot y$ are parallel,
and
$b \cdot u$ and $(bp^{-1}q) \cdot v$ are parallel.
Thus
$(a \cdot x)(b \cdot u)$
and
$((am^{-1}n) \cdot y)((bp^{-1}q) \cdot v)$
are defined and parallel.
We prove that
$((a \cdot x)(b \cdot u)) \equiv (((am^{-1}n) \cdot y)((bp^{-1}q) \cdot v))$.
Suppose that $\exists g \cdot ((a \cdot x)(b \cdot u))$.
Then by (A3), we have that
$\exists g \cdot (a \cdot x)$ and $g \cdot (b \cdot u)$.
By Lemma~2.1.1 it is easy to check that
$\exists g \cdot ((am^{-1}n) \cdot y)$
and
$\exists g \cdot ((bp^{-1}q) \cdot v)$.
It follows from (A2) that
$\exists g \cdot (((am^{-1}n) \cdot y)((bp^{-1}q) \cdot v))$.
To prove the converse,
suppose that
$\exists g \cdot (((am^{-1}n) \cdot y)((bp^{-1}q) \cdot v))$.
Then by (A3) we have that
$\exists g \cdot ((am^{-1}n) \cdot y)$
and
$\exists g \cdot ((bp^{-1}q) \cdot v))$.
From $\exists g \cdot ((am^{-1}n) \cdot y)$
we can deduce that $\exists (ga) \cdot x$ so that from $\exists a \cdot x$
and Lemma~2.1.1 we get $\exists g \cdot ( a \cdot x)$.
Similarly $\exists g \cdot (b \cdot u)$.
Thus by (A2) we have that $\exists g \cdot ((a \cdot x)(b \cdot u))$.
It follows that we have proved $[x][u] \E [y][v]$, as required.
Thus $\E$ is a congruence.

Finally, we prove that $\E$ is idempotent pure.
Suppose that $[x] \E [x][x]$.
By assumption $[x][x] = [(a \cdot x)(b \cdot x)]$.
By assumption for some $c$ and $d$ we have that
$c \cdot x \equiv d \cdot ((a \cdot x)(b \cdot x))$.
Thus $c \cdot x$ and $((da) \cdot x)((db) \cdot x)$ are
parallel.
Thus $\s (c \cdot x) = \s ((db) \cdot x)$.
Hence $c \cdot \s (x) = (db) \cdot \s (x)$.
But $G$ acts freely on $S_{o}$ and so $c = db$.
Similarly $c = da$.
Thus $da = db$ and so $a = b$.
It follows that $\s (a \cdot x) = \f (a \cdot x)$.
Hence $a \cdot \s (x) = a \cdot \f (x)$.
Hence $\s (x) = \f (x)$.
It follows that $x$ is a loop.
But by assumption $S$ is locally idempotent and so $x$ is idempotent,
which gives $[x]$ idempotent, as required.
It follows that $\E$ is idempotent pure.\qed

We call $\E$ the {\em empire congruence} on $S/G$.

\subsection{Universal groups}\label{sect2.2}
We now describe a procedure for
associating a group with a set
equipped with a partial binary operation $(S,\circ)$.
Let $FG(S)$ be the free group on the set $S$.
Elements of $FG(S)$ can be regarded as reduced strings
over the alphabet $S \cup S^{-1}$.
We shall represent such strings as lists.
There is therefore a function $\iota \colon \: S \rightarrow FG(S)$
given by $s \mapsto (s)$.
We denote the product in $FG(S)$ by $\cdot$.
Define the relation $\sim$ on $FG(S)$ to be the congruence
generated by all ordered pairs 
$((s) \cdot (t), (s \circ t))$
where $\exists s \circ t$.
Denote the $\sim$-congruence class containing $x \in FG(S)$ by $[x]$.
Let $\nu \colon \: FG(S) \rightarrow FG(S)/ \sim$
be the associated natural homomorphism.
Put $G(S) = FG(S)/\sim$ 
and $\tau = \nu \iota$.
It is easy to check that
$\tau \colon \: S \rightarrow G(S)$ is a morphism.
Furthermore,
if $\alpha \colon \: S \rightarrow G$ is any
morphism from $S$ to a group $G$,
then there is a unique homomorphism
$\theta \colon \: G(S) \rightarrow G$ such that
$\tau \theta = \alpha$.
It is therefore legitimate to call
$G(S)$ the {\em universal group} of $(S,\circ)$.\\

\noindent
{\bf Remark} Observe that the image
of $\tau$ generates $G(S)$.\\

In this paper, we are interested in
computing the universal groups
of some inverse semigroups associated with tilings
and point-sets.
The following results will be useful.

\begin{proposition}\label{prop2.2.1} Let $S$ and $T$ be inverse semigroups,
and let $\theta \colon \: S \rightarrow T$
be a surjective, idempotent pure homomorphism.
Then $S$ and $T$ have isomorphic universal groups.
\end{proposition}
\bew\ Let the universal morphisms be
$\alpha \colon \: S \rightarrow G(S)$ and
$\beta \colon \: T \rightarrow G(T)$.
By universal properties there is
a homomorphism $\theta^{\ast} \colon \: G(S) \rightarrow G(T)$
such that $\theta^{\ast}\alpha = \beta \theta$.

Define $\phi \colon \: T \rightarrow G(S)$
by $\phi (t) = \alpha (s)$ where $\theta (s) = t$.
Suppose that $\theta (s) = \theta (s')$.
Then $\theta (s)^{-1}\theta (s) = \theta (s)^{-1}\theta (s')$.
Hence $\theta (s^{-1}s) = \theta (s^{-1}s')$.
But $\theta$ is an idempotent pure homomorphism
and so $s^{-1}s'$ is defined and is an idempotent.
Thus $\alpha (s^{-1}s') = 1$ which gives
$\alpha (s) = \alpha (s')$.
It follows that $\phi$ is well-defined.

To show that $\phi$ is a morphism.
Let $t,t' \in T$ such that $tt'$ is defined.
Let $\theta (s) = t$ and $\theta (s') = t'$.
Because $\theta$ is a homomorphism $ss'$ is defined.
Hence result.

It follows by universal properties that there
is a unique homomorphism
$\phi^{\ast} \colon \: G(T) \rightarrow G(S)$
such that $\phi^{\ast}\beta = \phi$.

Let $\alpha (s) \in \mbox{im} \, \alpha$.
Then
$(\phi^{\ast}\theta^{\ast})(\alpha (s))
=
\phi^{\ast}\beta \theta (s)
=
\phi \theta (s)
=
\alpha (s)$.

Let $\beta (t) \in \mbox{im} \, \beta$.
Then
$(\theta^{\ast} \phi^{\ast})(\beta (t))
=
\theta^{\ast} \phi (t)
=
\theta^{\ast}\phi (\theta (s))
=
\theta^{\ast}\alpha (s)
=
\beta \theta (s)
=
\beta (t)$.
Thus
$\phi^{\ast}\theta^{\ast}$ is the identity
on the image of $\alpha$
and
$\theta^{\ast}\phi^{\ast}$ is the identity
on the image of $\beta$.
But the image of $\alpha$ generates $G(S)$
and the image of $\beta$ generates $G(T)$.
It follows that $\theta^{\ast}$ and $\phi^{\ast}$
are mutually inverse homomorphisms.
Thus $G(S)$ is isomorphic to $G(T)$,
as required.\qed

We now describe one method for constructing 
homomorphisms between strongly $E^{\ast}$-unitary semigroups.
Let $C$ and $D$ be inverse categories and let $G$ be a group 
which acts partially and freely on both $C$ and $D$.
A functor $\theta \colon \: C \rightarrow D$ is said to be
{\em $G$-equivariant} if the following condition holds:
if $x \in C$ then 
$\exists g \cdot x \Leftrightarrow \exists g \cdot \theta (x)$
and $\theta (g \cdot x) = g \cdot \theta (x)$.

\begin{proposition} Let $C$ and $D$ be
inverse categories on which $G$ acts partially and without fixed points.
Let $\theta \colon \: C \rightarrow D$ be a $G$-equivariant
functor. 
Define $[\theta] \colon \: C/G \rightarrow D/G$ by
$[\theta][x] = [\theta (x)]$.
Then $[\theta]$ is a homomorphism.
\end{proposition}
\bew\ Define $[\theta] \colon \: C/G \rightarrow D/G$ by
$[\theta][x] = [\theta (x)]$.
We claim this is a well-defined function:
for suppose $\exists g \cdot x$ then 
$\theta (g \cdot x) = g \cdot \theta (x)$ by the fact that
$\theta$ is $G$-equivariant.
Now $[x] = [g \cdot x]$
and $[\theta (g \cdot x)] = [g \cdot \theta (x)] = [\theta (x)]$.
Now we prove that $\theta$ is a homomorphism.
Suppose that $[x][y]$ is defined.
Then $[x][y] = [(a \cdot x)(b \cdot y)]$.
Thus $[\theta]([x][y]) = [\theta((a \cdot x)(b \cdot y))]
= [\theta (a \cdot x) \theta (b \cdot y)]$
since $\theta$ is a functor.
But then $\theta (a \cdot x) = a \cdot \theta (a)$
and $\theta (b \cdot y) = b \cdot \theta (y)$.
It follows that $[\theta (x)][\theta (y)]$ is defined
and we have proved that 
$[\theta]([x][y]) = [\theta]([x])[\theta]([y])$.
\qed

The most important class of inverse semigroups in this article
are the $F^{\ast}$-inverse semigroups.
Let $S$ be an $F^{\ast}$-inverse semigroup.
Let $M = M(S)$ be the set of maximal elements of $S$.
Define $\circ$ on $M$ as follows:
if $x,y \in M$ and $\exists xy$ then $xy$ is beneath
a unique maximal element $z$.
Define in this case $x \circ y = z$.

\begin{proposition}
Let $S$ be an $F^{\ast}$-inverse semigroup.
\begin{description}
\item[{\rm (i)}] The universal group of $S$ is the
same as the universal group of $(M(S),\circ)$.

\item[{\rm (ii)}] If $S$ is a monoid then $(M(S),\circ)$ is a group-like set.
\end{description}
\end{proposition}
\bew\ (i) Let $\alpha \colon \: S \rightarrow G$ be a morphism to a group.
We show first that $\alpha$ induces a morphism from $M(S)$ to $G$.
Suppose that $x \circ y$ is defined in $M(S)$.
Then $xy \leq x \circ y$ and so
$\alpha (xy) = \alpha (x \circ y)$.
But $xy$ defined means $\alpha (x)\alpha (y)$ defined
and so $\alpha (x \circ y) = \alpha (x) \alpha (y)$. 

Now suppose $\beta \colon \: M(S) \rightarrow G$ is a morphism.
We show that $\beta$ can be extended uniquely to $S$ to produce
a morphism from $S$ to $G$.
Let $a \in S$.
Then $a \leq x \in M(S)$ and $x$ is unique.
Define $\beta^{\ast}(a) = x$.
Clearly $\beta^{\ast}$ extends $\beta$.
Suppose $a,b \in S$ and $ab$ is defined.
Let $a \leq x$ and $b \leq y$.
Then $ab \leq xy \leq x \circ y$.
Thus $\beta^{\ast}(ab) = \beta (x \circ y)$.
Now $\beta^{\ast}(a) = \beta (x)$ and
$\beta^{\ast} (b) = \beta (y)$.
By assumption
$\beta (x)\beta (y) = \beta (x \circ y)$.
Hence
$\beta^{\ast}(ab) = \beta^{\ast}(a) \beta^{\ast}(b)$.
It follows that $\beta^{\ast}$ is a morphism
from $S$ to $G$ which extends $\beta$.
Suppose that $\gamma$ is a morphism from
$S$ to $G$ extending $\beta$.
Let $a \in S$ where $a \leq x \in M(S)$.
Then $\gamma (a) = \gamma (x) = \beta (x) = \beta^{\ast} (a)$.
Thus $\beta^{\ast}$ is unique.

We have essentially constructed a bijection
between the set of morphisms from $M(S)$ to $G$
and
the set of morphisms from $S$ to $G$.
It is now clear that the universal group of
$S$ is the same as the universal group of $M(G)$.

(ii) The proof of this assertion is straightforward
and is omitted.\qed

\subsection{Macbeath's theorem}
Let $G$ be a group acting partially on the set $V$.
We construct a group-like set from this partial action.
For each $g \in G$
such that $V \cap g^{-1}V \neq \emptyset$ define
$$\lambda_{g} \colon \:  
V \cap g^{-1}V
\rightarrow
V \cap gV$$
by
$\lambda_{g}(v) = gv$.
Put $M(G,V) = \{\lambda_{g} \colon \: \lambda_{g} \neq \emptyset\}$.
Define $\circ$ on $M(G,V)$ by 
$$\lambda_{g} \circ \lambda_{h} = \lambda_{gh}$$
iff $\lambda_{g}\lambda_{h} \neq \emptyset$.
The proof of the following is straightforward.

\begin{proposition}
$(M(G,V),\circ)$ is a group-like set.\qed
\end{proposition}

Observe that $\lambda_{g} \circ \lambda_{h}$ is defined
if and only if
$V \cap g^{-1}V \cap hV \neq \emptyset$
if and only if
$V \cap gV \cap ghV \neq \emptyset$.
Define
$$E(V) = \{g \in G \colon \: V \cap g^{-1}V \neq \emptyset \}$$
and define
$F'(V) \subseteq E(V) \times E(V)$ by
$$F'(V) = \{(g,g') \in G \times G \colon \: V \cap gV \cap gg'V 
\neq \emptyset \}.$$
Finally define $R(V)$ to be the set of ordered pairs
$$R(V) = \{(g \cdot g',gg') \colon \: (g,g') \in F' \}.$$
The pair $(E(V),R(V))$ defines a group $\mathbf{G}(G,V)$
which is the free group on $E(V)$ factored out
by the congruence $\equiv_{R(V)}$ generated by $R(V)$.
We say that the group $\mathbf{G}(G,V)$
is associated with the partial action of $G$ on $V$.
The following is now clear

\begin{proposition} Let $G$ act partially on the set $V$.
Then the universal group of the group-like set
$M(G,V)$ is isomorphic to $\mathbf{G}(G,V)$.\qed
\end{proposition}

There is one case where we can easily write down
the group associated with an action.

\begin{lemma}
Let $G$ act (globally) on the set $M$.
Then $E(M) = G$ and $F'(M) = G \times G$.
Thus in this case $\mathbf{G}(G,M)$ is isomorphic to $G$.
\end{lemma}
\bew\ The result is immediate because
$E(M) = G$ and $F'(M) = G \times G$ and so the corresponding presentation
is essentially just
the Cayley table of the group $G$.\qed

If $V \subseteq M$ is any subset then the action of
$G$ on $M$ induces a partial action of $G$ on $V$.

\begin{lemma}\label{lem2.3.4} Let $G$ act globally on $M$.
Let $U \subseteq V \subseteq W \subseteq M$ where $U$, $V$ and $W$ are
each equipped with the induced partial action.
\begin{description}

\item[{\rm (i)}] $E(U) \subseteq E(V) \subseteq E(W)$ 
and $F'(U) \subseteq F'(V) \subseteq F'(W)$. 
Hence $R(U) \subseteq R(V) \subseteq R(W)$.

\item[{\rm (ii)}] There are homomorphims
$\psi^{V}_{U} \colon \: \mathbf{G}(G,U) \rightarrow \mathbf{G}(G,V)$,
$\psi^{W}_{V} \colon \: \mathbf{G}(G,V) \rightarrow \mathbf{G}(G,W)$,
and
$\psi^{W}_{U} \colon \: \mathbf{G}(G,U) \rightarrow \mathbf{G}(G,W)$
such that
$\psi^{W}_{U} = \psi^{W}_{V}\psi^{V}_{U}$.

\end{description}
\end{lemma}
\bew\ (i) We prove that if $U \subseteq V$ then
$E(U) \subseteq E(V)$ and $F'(U) \subseteq F'(V)$.
Let $g \in E(U)$. 
Then $U \cap gU$ is nonempty.
But $U \cap gU \subseteq V \cap gV$ so that
$V \cap gV$ is non-empty.
The proof that $F'(U) \subseteq F'(V)$ is similar.
The fact that $R(U) \subseteq R(V)$ is now immediate.

(ii) Let $U \subseteq V$. Then by (i), 
we have that $E(U) \subseteq E(V)$ and $R(U) \subseteq R(V)$.
The function 
$\psi^{V}_{U} \colon \: \mathbf{G}(G,U) \rightarrow \mathbf{G}(G,V)$
is defined by $[g]_{R(U)} \mapsto [g]_{R(V)}$
where $[g]_{R}$ denotes the congruence class generated by $R$
containing $g$.
This is a well-defined homomorphism.
The proof of the remaining assertion is straightforward.\qed

\begin{lemma} Let $G$ act globally on the set $M$.
Let $V \subseteq M$ be equipped with the induced partial action
of $G$.
Suppose that $E(V)$ generates $G$. 
Then the map 
$\psi_V^M:\bG(G,V)\to \bG(G,M)$ is surjective.
\end{lemma}
\bew\  Let $[g]_{M} \in \bG(G,M)$.
By assumption, 
we can write $g = g_1\cdots g_n$ for some $g_i\in E(V)$. 
The product 
$[g_{1}]_{R(V)} \ldots [g_{n}]_{R(V)}$ is
a well-defined element of $\bG(G,V)$.
Its image under $\psi_{V}^{M}$ is $[g]_{M}$,
as required.\qed

Let $G$ be a topological group,
let $M$ be a hausdorff space,
and let $G \times M \rightarrow M$ be a continuous action
of $G$ on $M$.
Then $(G,M)$ is called a topological transformation group.
A subset $V \subseteq M$ is called a {\em $G$-covering}
if $M = \bigcup_{g \in G} gV$.
The following is Macbeath's Theorem~1 \cite{Macbeath64}.

\begin{theorem} Let $(G,M)$ be a topological transformation group
where $M$ is connected and simply-connected.
Let $V \subseteq M$ be open, path-connected and a $G$-covering.  
Then the homomorphism 
$\psi_V^M:\bG(G,V)\to \bG(G,M)$ 
is an isomorphism. \qed
\end{theorem}

The theorem above tells us that
the universal group of $M(G,V)$ is $G$ itself
under the given conditions.

\subsection{The universal group of $\Gamma (X,G,H)$}\label{sect1.4}
In this section, we shall obtain a necessary condition
for the universal group of $\Gamma = \Gamma (X,G,H)$
to be $G$ (Theorem~2.4.7).

Throughout this section we put $X' = X \cap G$ and recall that
we assume $1 \in X$.
By Theorem~2.1.2  the semigroups of the form $\Gamma$
are strongly $E^{\ast}$-unitary.
However a stronger result holds.

\begin{proposition} In the semigroup $\Gamma (X,G,H)$
the maximal elements are of the form
$[a,aX \cap bX,b]$ 
where $1 \in aX' \cap bX'$.
In particular, $\Gamma$ is an
$F^{\ast}$-inverse monoid.
\end{proposition}
\bew\  It is easy to check that
$[a,P,b] \leq [c,Q,d]$ if and only if
there exists $g \in G$ such that
$gc = a$, $gd = b$ and $P \subseteq gQ$.

Observe that if $[a,P,b]$ is an
element of $\Gamma$ then 
$1 \in P \subseteq aX \cap bX \subseteq aX,bX$
so that $1 \in aX \cap bX$;
this is equivalent to $a^{-1},b^{-1} \in X$,
but $a^{-1},b^{-1} \in G$ and so,
in turn, this is equivalent to $1 \in aX' \cap bX'$.
Thus $[a, aX \cap bX, b] \in \Gamma$ and
$[a,P,b] \leq [a,aX \cap bX, b]$.
However, the element $[a,aX \cap bX,b]$ is clearly maximal
using our characterisation of the natural partial order above.
It follows that the maximal elements of $\Gamma$ are the elements
of the form $[a,aX \cap bX,b]$ and
every element of $\Gamma$ is beneath a maximal element.
It is easy to check that each element of $\Gamma$ is beneath a unique 
maximal element.

The maximal idempotents have the form $[a, aX, a]$. 
However, any two such elements are equal.
It follows that $\Gamma$ is a monoid.\qed

By Proposition~2.2.3, it follows that the universal group of
$\Gamma (X,G,H)$ is the same as the universal group of the group-like
set of maximal elements of $\Gamma (X,G,H)$.
Thus to compute the universal group of $\Gamma$
we need to know more about these maximal elements.

\begin{lemma} Let $\Gamma = \Gamma (X,G,H)$.
Define $\phi \colon \: \Gamma \rightarrow G$ 
by $\phi [a,P,b] = a^{-1}b$.

\begin{description}

\item[{\rm (i)}] The function $\phi$ is a well-defined idempotent
pure morphism.

\item[{\rm (ii)}] The element $g \in \mbox{\rm im}(\phi)$
if and only if there exists $h \in X'$ such that $g^{-1}h \in X'$, 
in which case, 
$$\phi [h^{-1},h^{-1}X \cap h^{-1}gX,h^{-1}g] = g.$$

\item[{\rm (iii)}] $\mbox{\rm im}(\phi) = X'(X')^{-1}$.

\item[{\rm (iv)}] The restriction of $\phi$ to the set of
maximal elements of $\Gamma$ is injective.

\end{description}
\end{lemma}
\bew\ (i) This is straightforward.

(ii) To characterise the elements in $\mbox{im}(\phi)$,
it is enough, by Proposition~2.2.3, 
to find the elements of $G$ which are
the images of maximal elements in $\Gamma$ under $\phi$.
Thus $g \in \mbox{im}(\phi)$ if and only if
there exists a maximal element $[a,aX \cap bX,b]$ such that
$a^{-1}b = g$.
Thus  $g \in \mbox{im}(\phi)$ if and only if
$[a,aX \cap agX,ag]$ is maximal for some $a \in G$.
By Proposition~2.4.1,
$g \in \mbox{im}(\phi)$ if and only if
there exists $a \in G$ such that
$1 \in aX' \cap agX'$.
Thus $g \in \mbox{im}(\phi)$ if and only if
there exists $a \in G$ such that
$a^{-1},g^{-1}a^{-1} \in X'$.
Put $h = a^{-1}$, and we get the required result.

(iii) Suppose that
$g \in \mbox{im}(\phi)$.
Then by (ii) this implies there exists
$h \in X'$ such that $g^{-1}h \in X'$.
Put $g^{-1}h = k \in X'$.
Then $g = hk^{-1}$ where $h,k \in X'$.
Conversely, suppose that
$g = cd^{-1}$ where $c,d \in X'$.
Then $g^{-1}c = d$,
and so $c,g^{-1}c \in X'$
which by (ii) implies that $g \in \mbox{im}(\phi)$. 

(iv) Let $g \in \mbox{im}(\phi)$ and suppose that
$u,g^{-1}u \in X'$ and $v,g^{-1}v \in X'$ for some $u,v \in X'$.
Then by (ii), we have that
$$\phi [u^{-1}, u^{-1}X \cap u^{-1}gX,u^{-1}g]
= g =
\phi [v^{-1},v^{-1}X \cap v^{-1}gX,v^{-1}g].$$
However
$[u^{-1}, u^{-1}X \cap u^{-1}gX,u^{-1}g] 
= [v^{-1},v^{-1}X \cap v^{-1}gX,v^{-1}g]$.
It follows that each element of
$\mbox{im}(\phi)$ is the image
of {\em exactly one} maximal element of $\Gamma$.\qed

We can now obtain a more convenient description of the
universal group of the semigroup $\Gamma (X,G,H)$.

\begin{theorem} Let $\Gamma = \Gamma (X,G,H)$.
Then the universal group of $\Gamma$ is the 
group associated with the partial action of $G$ on $X'$.
\end{theorem}
\bew\ We prove that the group-like set  
$M(\Gamma)$ of maximal elements of $\Gamma (X,G,H)$
is isomorphic to the group-like set $M(G,X')$ (from Section~2.3).
This implies that their universal groups are isomorphic.
Hence the universal group of $\Gamma (X,G,H)$ is isomorphic
to the group ${\bf G}(G,X')$ (using Proposition~2.3.2).

Define 
$$\alpha \colon \: M(\Gamma) \rightarrow M(G,X')
\mbox{ by }
\alpha [a,aX \cap bX,b] = \lambda_{a^{-1}b}.$$
We show first that $\alpha$ is well-defined.
The element $[a,aX \cap bX,b]$ is maximal,
so that $1 \in aX' \cap bX'$ by Proposition~2.4.1.
Hence 
$$(a^{-1}b)^{-1}X' \cap X' = b^{-1}aX' \cap X' \neq \emptyset.$$
It follows that $\lambda_{a^{-1}b}$ is non-empty.
In addition, 
$a^{-1}b = \phi [a,aX \cap bX,b]$ is well-defined
by Lemma~2.4.2(i).
It follows that $\alpha$ is well-defined.

To show that $\alpha$ is injective,
suppose that
$$\alpha [a,aX \cap bX,b] = \alpha [c,cX \cap dX,d].$$
Then $\lambda_{a^{-1}b} = \lambda_{c^{-1}d}$.
However the partial action of $G$ on $X'$ is a restriction
of the left action of $H$ on itself.
It follows that $a^{-1}b = c^{-1}d$.
Injectivity now follows from the fact that $\phi$ is injective
when restricted to maximal elements by Lemma~2.4.2(iv).

To show that $\alpha$ is surjective,
let $\lambda_{g} \in M(G,X')$.
Then $g^{-1}X' \cap X' \neq \emptyset$.
It follows that there exists $h \in X'$ such that
$g^{-1}h \in X'$.
But then by Lemma~2.4.2,
we have that $g \in \mbox{\rm im}(\phi)$
and
$[h^{-1},h^{-1}X \cap h^{-1}gX,h^{-1}g]$
is a maximal element such that
$\phi [h^{-1},h^{-1}X \cap h^{-1}gX,h^{-1}g] = g$.
It follows that $\alpha$ is surjective.

We have therefore established a bijection between
$M(\Gamma)$ and $M(G,X')$.
We now show that this bijection underlies
an isomorphism of group-like sets.

Let $[a,aX \cap bX,b]$ and $[c,cX \cap dX,d]$
be elements of $M(\Gamma)$.
Let 
$$\phi [a,aX \cap bX,b] = a^{-1}b = g
\mbox{ and }  
\phi [c,cX \cap dX,d] = c^{-1}d = h.$$
We prove first that
$$\exists [a,aX \cap bX,b][c,cX \cap dX,d]
\Leftrightarrow
\lambda_{g}\lambda_{h} \neq \emptyset.$$ 
Suppose that
$[a,aX \cap bX,b][c,cX \cap dX,d]$ is defined.
Then it is easy to check that
$X' \cap gX' \cap ghX'$ is non-empty and
so $\lambda_{g}\lambda_{h}$ is non-empty.
Conversely, suppose that
$\lambda_{g}\lambda_{h}$ is non-empty.
Then there exists $p \in X' \cap gX' \cap ghX'$.
Define $n = p^{-1}a^{-1}bc^{-1}$ and $m = p^{-1}a^{-1}$.
Then $mb = nc$ and it is easy to check that
$$maX \cap mbX \cap ncX \cap ndX$$
contains 1.
It follows that
$$[a,aX \cap bX,b][c,cX \cap dX,d] 
=
[ma, maX \cap mbX \cap ncX \cap ndX,nd]$$
is defined.

Finally suppose that
$[a,aX \cap bX,b][c,cX \cap dX,d]$ is defined.
Then from the above
$$[a,aX \cap bX,b] \circ [c,cX \cap dX,d]
=
[ma, maX \cap ndX, nd]$$
where $mb = nc$.
Now
$\phi [ma, maX \cap ndX, nd] = (ma)^{-1}nd
= a^{-1}m^{-1}nd = a^{-1}(bc^{-1})d = gh$.
Hence 
$$\alpha ([a,aX \cap bX,b] \circ [c,cX \cap dX,d])
= \alpha [a,aX \cap bX,b] \circ \alpha [c,cX \cap dX,d].$$

We have therefore proved that $\alpha$ is a bijective morphism
of group-like sets with the extra property that
$$[a,aX \cap bX,b] \circ [c,cX \cap dX,d] 
\mbox{ is defined if and only if }
\alpha [a,aX \cap bX,b] \circ \alpha [c,cX \cap dX,d]
\mbox{ is defined}.$$
Hence $\alpha$ is an isomorphism of group-like sets.\qed

We shall now derive some sufficient conditions on $(X,G,H)$ in order that
the universal group of $\Gamma (X,G,H)$ be $G$.
We make the following blanket assumptions:

\begin{itemize}
\item $H$ is a connected and simply connected topological group.

\item $G$ is a dense subgroup of $H$. 

\item There is an open, path-connected subset $V \subseteq X$ which is
also a $G$-covering.
\end{itemize}

Put $X' = X \cap G$ and $V' = V \cap G$.
Because the conditions of Macbeath's theorem hold we know
that $\bG (G,V)$ is isomorphic to $G$.
We have also proved that the universal group of $\Gamma (X,G,H)$ is
isomorphic $\bG (G,X')$.
We aim to find conditions on $X$ in order that
$\bG (G,X')$ is isomorphic to $G$.
The following commutative diagram summarises the maps we are
interested in:

$$\diagram
\bG (G,V) \rto^{\psi^{X}_{V}} 
& \bG (G,X)\\
\bG (G,V') \uto^{\psi^{V}_{V'}} \rto^{\psi^{X'}_{V'}} 
& \bG (G,X') \uto_{\psi^{X}_{X'}}
\enddiagram$$

\begin{proposition} With reference to the above diagram, we have the following:
\begin{description}
\item[{\rm (i)}] $\psi_V^X:\bG(G,V)\to \bG(G,X)$ is injective. 
\item[{\rm (ii)}] $\psi^{V}_{V'} \colon \: \bG(G,V') \rightarrow \bG(G,V)$ is
an isomorphism.
\item[{\rm (iii)}] $\psi_{V'}^{X'}:\bG(G,V')\to \bG(G,X')$ is injective.
\end{description}
\end{proposition}
\bew\ (i) From $V \subseteq X \subseteq H$ and Lemma~\ref{lem2.3.4},
we have $\psi^{H}_{V} = \psi^{H}_{X} \psi^{X}_{V}$.
By Macbeath's theorem, $\psi^{H}_{V}$ is an isomorphism.
Thus $\psi^{X}_{V}$ is injective.

(ii) We show that $E(V') = E(V)$ and $F'(V') = F'(V)$
which immediately implies the result.
We already know that $E(V') \subseteq E(V)$
and $F'(V') \subseteq F'(V)$.
Let $g \in E(V)$.
Then $V \cap gV$ is non-empty by definition.
But $V$ is open and so $gV$ is open.
Thus $V \cap gV$ is a non-empty open set.
Now $G$ is dense in $H$ and so 
$G \cap V \cap gV$ is also non-empty.
But
$G \cap V \cap gV = G \cap V \cap gG \cap gV
= (G \cap V) \cap g(G \cap V)$.
Thus $g \in E(V')$.

We already know that $F'(V') \subseteq F'(V)$.
We prove that $F'(V) \subseteq F'(V')$.
Let $(g,h) \in F'(V)$.
Then $V \cap gV \cap ghV$ is non-empty by definition.
It is also open.
Thus by denseness
$G \cap V \cap gV \cap ghV$ is non-empty.
It follows that $(g,h) \in F'(V')$.

(iii) Suppose that $G$ is a dense subgroup of $H$.
Since $V' \subseteq X'$ there is a homomorphism
$\psi^{X'}_{V'} \colon \: \bG(G,V') \rightarrow \bG (G,X')$.
Now $\psi_V^X \psi^{V}_{V'} = \psi_{X'}^{X} \psi_{V'}^{X'}$
by Lemma~\ref{lem2.3.4}.
Using (i) and (ii) of this lemma we see that
the lefthand side is injective.
It follows that  $\psi^{X'}_{V'}$ is injective.
\qed

We find a condition which guarantees that
$\psi_{V'}^{X'}:\bG(G,V')\to \bG(G,X')$ is surjective. 
We begin with a lemma.

\begin{lemma} Let $1 \in V$.
Let $a,b \in V$.
Then there exists an open set $U$ such that
$1 \in U \subseteq V$ such that
$Ua \cup Ub \subseteq V$.
\end{lemma}
\bew\ Put $U = (Va^{-1} \cap Vb^{-1}) \cap V$.
Since $1,a,b \in V$ we have that $1 \in U$.
Since $V$ is open it follows that $Va^{-1}$ and $Vb^{-1}$ are open.
Thus $U$ is an open set satisfying the conditions.\qed

\begin{proposition} Suppose that
$1\in V\subseteq X \subseteq \bigcup_{s\in V} Vs^{-1}$. 
Then 
\begin{description}
\item[{\rm (i)}] $\psi^{X'}_{V'}$ is surjective.
\item[{\rm (ii)}] $\psi_V^X$ is surjective
\end{description}
\end{proposition}
\bew\ (i) We have that $1 \in V' \subseteq X' \subseteq X$.
Let $g \in E(X')$.
Then by definition $X' \cap gX' \neq \emptyset$.
Thus there exists $x \in X'$ such that $g^{-1}x \in X'$.
By assumption, there exist $s,t \in V$
such that $x \in Vs^{-1}$ and $g^{-1}x \in Vt^{-1}$.
Thus $V \cap x^{-1}V$ is non-empty and an open set.
Thus $V \cap x^{-1}V \cap G$ is non-empty.
It follows that $x \in E(V')$.
Similarly $g^{-1}x \in E(V')$.
Finally, $(x,x^{-1}g) \in F'(X)$
because $x, x^{-1}x, g^{-1}x \in X'$.

(ii) Let $g \in E(X)$.
Then by definition $X \cap gX \neq \emptyset$.
Thus there exists $x \in X$ such that $g^{-1}x \in X$.
By assumption, there exist $s,t \in V$
such that $x \in Vs^{-1}$ and $g^{-1}x \in Vt^{-1}$.
By Lemma~2.4.5, we can find an open set $U$ such that
$1 \in U \subseteq V$ and $Us \cup Ut \subseteq V$.
It follows that for each $\epsilon \in U$ we have that
$\epsilon s, \epsilon t \in V$.
Now $U$ open implies $U^{-1}$ open and so $xU^{-1}$ is open.
But $G$ is dense and so $xU^{-1} \cap G$ is non-empty.
Choose $\epsilon \in U$ such that
$x\epsilon^{-1} \in G$.
Now $xs \in V$ and $\epsilon s \in V$.
Thus $s \in x^{-1}V$ and $s \in \epsilon^{-1}V$.
Hence $V \cap x \epsilon^{-1}V \neq \emptyset$.
Thus $x \epsilon^{-1} \in E(V)$.
Similarly $g^{-1}xt, \epsilon t \in V$ and so
$t \in (g^{-1}x)^{-1}V \cap \epsilon^{-1} V$.
Hence $\epsilon x^{-1}g V \cap V \neq \emptyset$.
Thus $\epsilon x^{-1} g \in E(V)$.
It follows that $[x\epsilon^{-1}]_{R(V)}$ and
$[\epsilon x^{-1}g]_{R(V)}$ are well-defined elements
of $\bG (G,V)$.
Observe that $x, \epsilon , g^{-1}x \in X$, 
Thus $x \in X \cap x \epsilon^{-1}X \cap gX$.
Hence $(x\epsilon^{-1}, \epsilon x^{-1}g) \in R(X)$.
Thus
$$\psi^{X}_{V}([x\epsilon^{-1}]_{R(V)}[\epsilon x^{-1}g]_{R(V)})
=
[x\epsilon^{-1}]_{R(X)}[\epsilon x^{-1}g]_{R(X)}
=
[g]_{R(X)}$$
as required. \qed

The proof of the following now follows from 
Propositions~2.4.4 and 2.4.6 combined with Macbeath's Theorem.

\begin{theorem} Let $H$ be a connected and simply connected topological group
and let $G$  be a dense subgroup of $H$. 
Let $1 \in V \subseteq X$ 
where $V$ is open in $H$, 
path-connected and a $G$-covering.
Suppose that $X \subseteq \bigcup_{s\in V} Vs^{-1}$. 
Then the universal group of $\Gamma(X,G,H)$ is isomorphic to $G$.\qed
\end{theorem}

The following example shows that the theorem fails to be true if we
merely require that $X$ contains an open path-connected set.
\bigskip

\noindent 
{\bf Example}
Let $H=\RR^n$, $G$ a dense subgroup of $H$ and $V$ an open
subset of $\RR^n$ of diameter less than $r$ and $X=V\cup (V+s)$ where
$s\in G$ with $|s|>8r$. Note that $s\notin E(V)$.
We claim that $\psi_V^{X'}$ is not surjective. 

To prove the claim we construct
a map $\chi:E(X')\to \Z$ as follows: 
if $x$ is a real number which is not in $\Z+\frac{1}{2}$
let $\{x\}$ be the integer number which is closest to $x$. Let
$q:\RR^n\to \RR$ be the orthoprojection onto the linear sub-space spanned by
$s$ and define for $x\in E(X')$
$$\chi(x) = \left\{\frac{q(x)}{|s|}\right\}.$$   
This works since $\frac{q(E(X'))}{|s|}\subseteq B_{\frac{1}{4}}(-1)\cup
B_{\frac{1}{4}}(0) \cup B_{\frac{1}{4}}(1)$ where $B_{\frac{1}{4}}(t)$
is the open ball of radius $\frac{1}{4}$ around $t\in \RR$. Now it is clear
that for $x,y,xy\in E(X')$ holds

$$\left\{ \frac{q(xy)}{|s|}\right\}=\left\{\frac{q(x)+q(y)}{|s|}\right\}=
\left\{\frac{q(x)}{|s|}\right\}+\left\{\frac{q(y)}{|s|} \right\}.$$

Hence if $(x,y)\in F'(X')$ then $\chi(xy)=\chi(x)+\chi(y)$ which implies
that $\chi$ induces a well-defined homomorphism $\bG(G,X')\to
\Z$. Since $s\in E(X')$ and $\chi(s)=1$ this
homomorphism is surjective.

Clearly $\chi$ maps $\im\psi_V^{X'}$ onto $\{0\}$ showing that $\psi_V^{X'}$
cannot be surjective. \qed

\section{Point-sets}\setcounter{theorem}{0}
In this section, we shall investigate more closely the algebraic
description of point-sets of $\RR^d$. We elaborate on the definition
of their associated semigroup and compute examples of universal groups
for such semigroups. The emphasis in the last subsection lies on model
sets which we show to be closely linked to the 
semigroups of Section~\ref{sect1.4}. 

\subsection{Point-set semigroups}
Any subset of $\mathbb{R}^{n}$ is called a {\em point-set};
it is set to be {\em discrete} 
if its intersection with any closed ball is a
finite set;
it is said to be {\em uniformly discrete} if there exists an
$r>0$ such that $|x - y| \geq r$
for all its points $x,y$; 
it is said to be {\em relatively dense}
if there exists an
$R>0$ such that every sphere of radius $R$  
contains at least one of its points.
A {\em Delone set} is a subset of $\mathbb{R}^{n}$
which is both 
uniformly discrete and relatively dense.
In physics such sets are used to describe the (equilibrium) positions
of atoms in an (infinitely extended) solid.

If $\La \subseteq \mathbb{R}^{n}$ is a point-set,
then we can construct the point-set semigroup $\Gamma (\La)$
according to the method of Example~(iii) following Theorem~2.1.2.
We recall the construction in more detail now.

Let $\La$ be a subset of $\mathbb{R}^{n}$. 
Put
$$C(\La) = \{(p_2,P,p_1)\colon \: 
p_i\in P\subseteq \La,
\mbox{ where }
P \mbox{ is finite}\}.$$ 
It is easy to check that
$C(\La)$ is an inverse category which is locally idempotent
when we define a partial product as follows
$$(p_{2},P,p_{1})\circ (q_{2},Q,q_{1}) = 
\left
\{\begin{array}{ll}
(p_{2},P \cup Q,q_{1}) & \mbox{if $p_{1} = q_{2}$}\\
\mbox{undefined}       & \mbox{else.}  
\end{array}
\right.
$$
Right and left identities are ${\bf d}(p_{2},P,p_{1}) = (p_{1},p_{1},p_{1})$
and ${\bf r}(p_{2},P,p_{1}) = (p_{2},p_{2},p_{2})$ and the inverse of
$(p_{2},P,p_{1})$ is $(p_{1},P,p_{2})$.

We may define a partial action of $\mathbb{R}^{n}$ on $C(\La)$ 
as follows:
$\exists g\cdot (p_2,P,p_1)$ 
whenever $P + g \in \La$ 
and then 
$$ g\cdot (p_2,P,p_1)= (p_2 + g,P + g,p_1 + g).$$
The set of equivalence classes under this action
is denoted by
$\Gamma (\La) = C(\La)/\mathbb{R}^{n}$,
and this is a strongly $E^{\ast}$-unitary inverse semigroup.
In fact, we can say a little more about these semigroups.

\begin{proposition}\label{prop3.1.1}
Point-set semigroups are strongly $F^{\ast}$-inverse monoids.
\end{proposition}
\bew\ Let $\La$ be a point-set and $\Gamma (\La)$ its
associated semigroup.

Observe first that the non-idempotent maximal elements of $\Gamma (\La)$
are those elements of the form $[y,\{y,x\},x]$,
since we have $[y,P,x] \leq [y,\{ y,x\},x]$.
There is exactly one maximal idempotent,
namely $[x,\{x \},x]$, since any two elements of $\D$
are related by a translation. 
Thus $[x,\{x \},x]$ is the identity.
It is now clear that $\Gamma (\La)$ is $F^{\ast}$-inverse.

Finally, we show explicitly 
that $\Gamma (\La)$ is strongly $E^{\ast}$-unitary. 
Let 
$$H_{\La} = \langle \La - \La \rangle,$$ 
the subgroup of $\mathbb{R}^{n}$ 
generated by the vectors of differences 
$\{p-q \colon \: p,q \in \La \}$.
The map 
$\varphi \colon \: \Gamma (\La) \rightarrow H_{\La}$
given by
\begin{equation}\label{eq1}
\varphi[p_{2},P,p_{1}]=p_2 - p_1
\end{equation}
is a surjective, idempotent pure morphism.\qed 

Point-set semigroups are $F^{\ast}$-inverse monoids.
Thus their universal groups are the universal groups of
their group-like sets of maximal elements by Proposition~2.2.3.
We shall now  obtain a direct description of the group-like
set of maximal elements of $\Gamma (\La)$.

Consider the set $\La - \La$.
This is furnished with a partial binary operation as follows.
Let $a,b \in \La - \La$.
Suppose that $a = x - y$ and $b = y - z$ for some
$x,y,z \in \La$.
Then define $a \oplus b = a + b$,
which is clearly an element of $\La - \La$.
It is evident that $\oplus$ is a well-defined operation
on $\La - \La$ and that $(\La - \La,\oplus)$
is a group-like set.

\begin{proposition}\label{prop1} Let $\La$ be a point-set.
Then the group-like set of maximal elements of $\Gamma (\La)$
is isomorphic to $\La - \La$ equipped with the operation $\oplus$.
\end{proposition}
\bew\ The maximal elements of $\Gamma (\La)$ have the form
$[y, \{y,x\},x]$ where $x,y \in \La$.
Denote the set of maximal elements of $\Gamma (\La)$ by $M(\La)$.
Define 
$$\theta \colon \: \La - \La \rightarrow M(\La)
\mbox{ by } 
\theta (y-x) = [y,\{y,x\},x].$$
We first show that $\theta$ is well-defined.
Suppose that $y - x = v - u$ where $x,y,u,v \in \La$.
Put $a = u - x$.
Then 
$$a \cdot (y, \{y,x \},x)
= (a + y, \{a + y, a + x \}, a + x)
= (v, \{v,u \},u).$$
Hence $[y, \{y,x\},x] = [v, \{v,u \},u]$, as required.
It is clear that $\theta$ is injective,
and immediate that $\theta$ is surjective.

Let $a,b \in \La - \La$ and suppose that
$a \oplus b$ is defined.
Then there exist $x,y,z \in \La$ such that
$a = x - y$ and $b = y - z$.
By definition, $a \oplus b = x - z$.
Now
$\theta (a) = [x, \{x,y \},y]$
and
$\theta (b) = [y, \{y,z\},z]$.
We have that 
$$\theta (a) \theta (b) = [x \{x,y,z\}, z]
\leq [x, \{x,z\}, z].$$
Hence $$\theta (a) \circ \theta (b) = [x ,\{x,z\},z] =
\theta (a \oplus b)$$
as required.
Thus $\exists a \oplus b$ implies $\exists \theta (a) \circ \theta (b)$
and $\theta (a \oplus b) = \theta (a) \circ \theta (b)$.
 
Now let $a, b \in \La - \La$ be such that
$\exists \theta (a) \circ \theta (b)$.
Let $a = x - y$ and $b = w - z$ where
$x,y,w,z \in \La$.
By assumption, 
$\exists [x, \{x,y\}, y][w,\{w,z \},z]$.
Thus there are translations
$g$ and $h$ such that
$g + y = h + w \in \La$
and $g + x, h + z \in \La$.
It follows that
$$a = x - y = (g + x) - (g + y)
\mbox{ and }
b = w - z = (h + w) - (h + z).$$
Thus $\exists a \oplus b$,
and $\theta (a \oplus b) = \theta (a) \circ \theta (b)$.

Finally, observe that $\theta (0) = [0,\{ 0\},0]$,
and 
$$\theta (-(y - x)) = \theta (x - y) = [x,\{x,y \},y]
= [y,\{x,y\},x]^{-1}.$$
Thus $\theta$ is an isomorphism of group-like sets.\qed 

There are now two groups naturally associated with a point-set $\La$: 

\begin{itemize}
\item The universal group of the point-set semigroup of $\La$
denoted by $G_{\La}$, which, by Proposition~3.1.2, is isomorphic
to the universal group of $(\La - \La,\oplus)$.

\item The group $H_{\La}$ which appeared in the proof of
Proposition~\ref{prop3.1.1}: it is the subgroup of $\mathbb{R}^{n}$
generated by $\La - \La$. 
\end{itemize}
Our aim is to study both $G_{\La}$ and its relation to $H_{\La}$.

\subsection{Lagarias' hierarchy of Delone sets}

Lagarias \cite{Lagarias} set up a hierarchy of Delone sets 
by imposing conditions on the set of difference vectors. 
This arose from an attempt to study the right conditions
on a Delone set required to describe quasicrystals. 
This forms part of a general programme in
understanding the mathematical conditions needed
for a set to describe an ordered solid.  
Loosely speaking, a Delone set is considered to be more ordered the
further down it is in the following list. 
A Delone set $\D$ is said to be:

\begin{description}
\item[{\rm (1)}] {\em finitely generated} if $H_\D$ is finitely generated,
\item[{\rm (2)}] of {\em finite local complexity}\footnote{In the literature one
    finds different terms for this property. Lagarias used {\em finite type}.} 
if $\D-\D$ is discrete,
\item[{\rm (3)}] a {\em Meyer set} if $\D-\D$ is a Delone set, 
\item[{\rm (4)}] a {\em model set} if it is defined by 
a cut-and-projection scheme (see below).
\end{description}

Delone sets with property (n) form a subclass of Delone sets with
property (n-1). 
The following theorem shows that the condition on $\D$
that $G_\D$ be finitely generated falls between (1) and (2).
We do not know, however, whether Delone sets with finitely generated
universal groups form a distinct class.

\begin{theorem}
If $\D$ is a Delone set of finite local complexity then $G_\D$ is
finitely generated. 
A Delone set for which $G_\D$ is
finitely generated is finitely generated.
\end{theorem}
\bew\ $\D$ is relatively dense so let $R>0$ such that any 
$R$-ball contains a point of $\D$ and set 
$$\Delta=\{\xi\in\D-\D \colon \: |\xi|\leq R\}.$$
Let $\xi=y-x$, $x,y\in \D$. Since any $R$-ball contains a point of
$\D$ we can find a sequence of points $x_0,\cdots,x_k\in\D $ such
that $x_0=x$, $x_k=y$ and $|x_i-x_{i-1}|\leq R$. Taking $\xi_i=
x_i-x_{i-1}$ we have $\xi_i\in\Delta$ and
$x+\sum_{i=1}^j \xi_i=x_j\in\D$ for all $j\leq k$. Hence
$(y-x)=(\xi_1\oplus\xi_2\cdots \oplus \xi_k)$.
Therefore, $G_\D$ is generated by the image of $\Delta$ in $G_\D$.
If $\D$ has finite local complexity then $\Delta$ is discrete and
hence a finite set which implies
that $G_\D$ finitely generated.

The second statement follows from the observation 
that the map in (\ref{eq1}) 
is a morphism so that by the universal property 
$H_\D$ is a homomorphic image of $G_\D$. \qed

\subsection{Point-sets in $\mathbb{R}$}\setcounter{theorem}{0}
We shall now investigate infinite discrete subsets of $\mathbb{R}$
with a view to carrying out some explicit calculations.
 
Let $\D \subseteq \mathbb{R}$ be a discrete subset. 
Choose $r_0\in \D$ and write $\D=\{r_i \colon \: i \in \Z\}$
where $i < j$ implies that $r_{i} < r_{j}$.
It follows that the set $\Sigma=\{r_{i}- r_{i-1} \colon \: i\in \Z\}$
is a set of distinct positive real numbers which is contained in $\D-\D$.
Each positive element of $\D-\D$ can be written as a sum of elements
from $\Sigma$.
It follows that
$$\Sigma \subseteq \D - \D \subseteq \erz{\Sigma}.$$
This implies that $H_{\D}$ is the subgroup of $\mathbb{R}$ 
generated by $\Sigma$.

We are interested in calculating $G_{\D}$.
To do this it is useful to relate the point-set $\D$ to a suitable tiling.
Regard $\Sigma$ as an alphabet equipped with a function
$a \mapsto |a|$ which gives the {\em length of $a$}:
if $a \in \Sigma$ 
corresponds to 
$r_{i} - r_{i-1}$ then $|a| = r_{i} - r_{i-1}$. 
Observe that by assumption if $a,b \in \Sigma$ and $a \neq b$ then
$|a| \neq |b|$.
Define a bi-infinite string $\T \colon \: \Z \to \Sigma$
by $\T (i) = a \in \Sigma$ if $|a| = r_{i} - r_{i-1}$.

Conversely, the tiling $\T$ and the length function $||$
determine $\D$ {\em up to translation}:
$\T$ determines a discrete subset $\{r_i \colon \: i\in\Z\}\subseteq \RR$ 
if we fix a value for $r_0$ and then set $r_i-r_{i-1}$ to be the real number
corresponding to $\T (i)$.
Knowledge of $r_{0}$ would enable us to recapture the original
point-set $\D$ exactly.
However, this is not necessary, 
because two point-sets which differ by a translation must have the
same universal groups $G_{\D}$.

Given $\D$, we fix an alphabet $\Sigma$, 
tiling $\T$ and length function $||$.
We show how to compute $G_{D}$ from $(\T,||)$.

It is useful to extend $||$ to any string over $\Sigma$:
if $u=a_1\dots a_n$ is a non-empty string over $\Sigma$,
define $|u|=\sum_i |a_i|$ to be its {\em length};
the empty string is defined to have length zero.

\begin{proposition} The group $G_{\D}$ is isomorphic to the free group
on $\Sigma$ factored out by the congruence generated
by all pairs of strings $(u,v)$ such that
$u$ and $v$ occur in $\T$ and have the same length.
\end{proposition}

\bew\ By Proposition~3.1.2, the group $G_{\D}$ is the universal group
of the group-like set $(\D - \D, \oplus)$. 

The elements of $\D - \D$ are either positive, negative, or zero.
Let $P$ be the subset of $\D - \D$ consisting of 0 and the positive elements,
the `positive' elements,
let $P^{-}$ be the remaining elements of $\D - \D$,
the `negative' elements.
The set $P$ inherits the partial binary operation $\oplus$.
We claim that the universal group of $(P,\oplus)$ is the same
as the universal group of $(\D - \D,\oplus)$.
Let $a,b \in \D - \D$ such that $\exists a \oplus b$.
There are four types of such such sums:
\begin{description}
\item[{\rm (1)}] $a,b \in P$ and so $a \oplus b \in P$.
\item[{\rm (2)}] $a,b \in P^{-}$ and so $a \oplus b \in P^{-}$;
observe that $a \oplus b = -(-b \oplus -a)$.
\item[{\rm (3)}] $a \in P$ and $b \in P^{-}$: if $|a| \geq |b|$
(absolute values) then $a \oplus b \in P$
giving $a = (a \oplus b) \oplus (-b)$;
else $a \oplus b \in P^{-}$
giving $-b = -(a \oplus b) \oplus a$.
\item[{\rm (4)}] $a \in P^{-}$ and $b \in P$: if $|a| > |b|$
then $a \oplus b \in P^{-}$
giving $-a = (-a \oplus -b) \oplus b$;
else $a \oplus b \in P$
giving $b = (a \oplus b) \oplus -a$.
\end{description}
Let $\alpha$ be a morphism from $(P,\oplus)$ to a group $G$.
Because $\exists 0 \oplus 0 = 0$,
it is immediate that $\alpha (0) = 1$,
the identity of $G$.
We show that $\alpha$ can be uniquely extended to a morphism $\alpha^{\ast}$
from the group-like set $(\D - \D, \oplus)$ to $G$.
Define 
$$\alpha^{\ast}(a) = \left\{ \begin{array}{ll}
\alpha (a) &\mbox{if $a \in P$}\\
\alpha (-a)^{-1} &\mbox{if $a \in P^{-}$}
\end{array}
\right.
$$
Our claim will be vindicated if we can prove that
$\alpha^{\ast}$ is a morphism,
but this is straightforward to check using
the four forms the product can take above.
It follows that $G_{\D}$ is the same as the universal group
of $(P,\oplus)$.

We now obtain a more convenient description of $(P,\oplus)$.
Let $L \subseteq \Sigma^{\ast}$ consist of all strings which occur in $\T$
together with the empty string.
Define an equivalence relation on $L$ by deeming
two strings $x,y \in L$ to be equivalent
if and only if $|x| = |y|$.
Denote the equivalence class containing $x$ by $[x]$,
and denote the set of equivalence classes by $S$.
Define a partial binary operation $\circ$ on $S$ as follows:
$[x] \circ [y]$ is defined iff there exist
$u \in [x]$ and $v \in [y]$ such that $vu \in L$,
in which case $[x] \circ [y] = [vu]$.
It is easy to check that this operation is well-defined.
We claim that $(S,\circ)$ and $(P,\oplus)$ are isomorphic.
Define $\theta \colon \: S \rightarrow P$ by $[x] \mapsto |x|$.
This is clearly well-defined, injective
and, from the definition of $\T$, it is also surjective.
We show $\theta$ is a morphism.
Suppose $\exists [x] \circ [y]$.
Then $[x] \circ [y] = [vu]$
where $u \in [x]$, $v \in [y]$ and $u,v,vu \in L$.
From the definition of the tiling $\T$,
it follows that there exist $r_{i},r_{j},r_{k} \in \D$
such that 
$|u| = r_{i} - r_{j}$, $|v| = r_{j} - r_{k}$
and 
$u = \T (j + 1) \ldots \T (i)$,
$v = \T (k+1) \ldots \T (j)$. 
Now $vu = \T (k+1) \ldots \T (i)$ and
$|u| \oplus |v| = r_{i} - r_{k}$ exists.
Thus $\theta ([x]) \oplus \theta ([y])$ exists and
equals $|u| \oplus |v|$
whereas 
$$\theta ([x] \circ [y]) 
= \theta ([vu]) = |vu| 
= |v| + |u|
= |u| \oplus |v|,$$
as required.

To show that $\theta$ is a homomorphism
suppose that $\theta ([x]) \oplus \theta ([y])$ is defined.
Then there exist $r_{i},r_{j},r_{k} \in \D$
such that
$|x| = r_{i} - r_{j}$ and $|y| = r_{j} - r_{k}$.
Let 
$u = \T (j + 1) \ldots \T (i)$
and
$v = \T (k+1) \ldots \T (j)$. 
By assumption, $vu \in \T$
and $|u| = |x|$ and $|v| = |y|$.
It follows that $[x] \circ [y]$ is defined in $S$.

We have therefore proved that $\theta$ is a bijective
homomorphism.
Thus $(P,\oplus)$ is isomorphic to $(S,\circ)$.

We now explicitly compute the universal group of $(S,\circ)$.
Let $FG = FG(\Sigma)$ be the free group on $\Sigma$:
elements of $FG$ are represented by lists.
Thus $a \in \Sigma$ maps to $(a)$ in $FG$.
The product in $FG$ is denoted by $\cdot$.
Let $\lambda$ be the congruence on $FG$
generated by pairs $(u,v)$ such that $u,v \in L$ and $|u| = |v|$.
Put $G = FG/\lambda$.
Define the function $\alpha \colon \: S \rightarrow G$
by $[x] \mapsto \lambda (x)^{-1}$.
This is well-defined, 
for suppose $[x] = [y]$ then $|x| = |y|$ and $x,y \in L$;
thus $\lambda (x) = \lambda (y)$
and so $\lambda (x)^{-1} = \lambda (y)^{-1}$
giving $\alpha ([x]) = \alpha ([y])$.
The function $\alpha$ is a morphism,
for suppose $\exists [x][y]$.
Then there exists $u \in [x]$ and $v \in [y]$ such that
$vu \in L$.
Now $\alpha ([x]) = \lambda (u)^{-1}$
and
$\alpha ([y]) = \lambda (v)^{-1}$
and 
$$\alpha ([x] \circ [y]) 
=
\alpha ([vu])
=
\lambda (vu)^{-1}
= 
\lambda (u)^{-1} \lambda (v)^{-1}
= \alpha ([x]) \alpha ([y]).$$

It remains to check that $\alpha$ is universal.
Let $\beta \colon \: S \rightarrow H$ be any morphism to a group.
Observe that every non-empty element of $S$ can be written
as a product of elements of the form $[a]$ where $a \in \Sigma$.
Indeed if $x = a_{1} \ldots a_{n} \in L$,
then $[x] = [a_{n}] \circ \ldots \circ [a_{1}]$.
Define $\theta \colon \: \Sigma \rightarrow H$
by $\theta (a) = \beta ([a])^{-1}$.
Then because $FG$ is the free group on $\Sigma$,
we can extend $\theta$ to a homomorphism from $FG$ to $H$.
Let $x \in FG(\Sigma)$ where $x \in L$
(so we think of $u$ as a reduced string),
and let
$x = a_{1} \ldots a_{n}$.
Then 
$$\theta (x) 
= \theta (a_{1}) \ldots \theta (a_{n})
= \beta ([a_{1}])^{-1} \ldots \beta ([a_{n}])^{-1}
= (\beta ([a_{n}]) \ldots \beta ([a_{1}]))^{-1}
= \beta ([a_{1} \ldots a_{n}])^{-1}
= \beta ([x])^{-1}.$$
Let $u,v \in L$ such that $|u| = |v|$.
Then $[u] = [v]$ and so $\beta ([u]) = \beta ([v])$
but this implies
$\theta (u) = \theta (v)$.
Thus $\theta$ induces a homomorphism
$\gamma$ from $G$ to $H$.
By construction $\gamma \alpha = \beta$
(the two inverses cancel);
it is easy to see that $\gamma$ is the unique
homomorphism with this property.\qed

\noindent
{\bf Examples}\\

We now use the description of discrete subsets by bi-infinite
strings and the last proposition to provide a table of simple
examples, all being based on a two-letter alphabet $\Sigma=\{a,b\}$. 
In particular, in all cases $G_\D$ will be a homomorphic image of $FG_2$, 
the free group on two generators. 
Now $G_\D$ depends not only on the bi-infinite string $\mathcal T$ 
which represents it, but also on the lengths of $a$ and $b$. 
It is, however, only the ratio $\frac{|a|}{|b|}$ which is important:
the group $G_{\D}$ is unchanged if we define new length functions
which are simply multiples of the original by a fixed positive number.
In all cases, we shall assume that $\T$ contains an infinite number 
of $a$'s and $b$'s and that the ratio is not $1$ (this simply means
that $|a| \neq |b|$).

\begin{description}
\item[{\rm (Case~1).}] This is the generic case in which
the string representing $\D$ contains both $ab$ and $ba$. 
It follows that $(a)$ and $(b)$ commute and $G_\D$ is a quotient of $\Z^2$. 
Moreover, only if $\frac{|a|}{|b|}\in\Q$ can there be
two strings in $\T$ which have the same length 
but a different number of $a$'s or $b$'s.\\

\item[{\rm (Case~2).}] Here $\D$ is a periodic repetition of 
$ab$ and $\frac{|a|}{|b|} \in \Q$. 
We claim that there is no other relation between
$(a)$ and $(b)$ than their commutativity. 
To see why, note first that if a string from the tiling
has $n$ $a$'s and $m$ $b$'s then $ n = m \pm 1$.
Thus the equality of the lengths of two strings implies
that we have an equation $n|a|+m|b|=k|a|+l|b|$ 
where $n=m\pm 1$ and $k=l\pm 1$
and $n,m,k,l$ natural numbers ($|a|,|b|$ are positive). 
Since $|a| \neq |b|$ the only solution to this equation is $n=k, m=l$.\\

\item[{\rm (Case~3).}] Here $\D$ consists of a half infinite $a$-sequence
matched to a half infinite $b$-sequence. If $\frac{|a|}{|b|}$ is
irrational there are no distinct strings of equal length. 
This provides an example of a non-abelian universal group.
A more complicated group arises if $\frac{|a|}{|b|}=\frac{n}{m}$, 
$n,m\in\Nat$ coprime. In this case, only the pairs of strings 
$|a|^{mr+i}|b|^j$ and
$|a|^i|b|^{nr+j}$, $i,j\in\Nat$,  
have the same length.
Hence $(a)^m=(b)^n$ which then implies that this element
commutes with everything.
\end{description}

\begin{table}[ht]
\caption{Universal group and group of difference vectors for point-sets 
in $\RR$}\label
{eqtable}
\medskip
\renewcommand\arraystretch{1.5}
\begin{center}
\begin{tabular}{|c|c|c|l|l|}
\hline
& $\mathcal T$ &$\frac{|a|}{|b|}$ & $G_\D$ & $H_\D$ \\
\hline \hline
case 1& $ab, ba\in\mathcal T_\D$ &$\notin \Q$ &  $\Z^2$& $\Z^2$\\ 
\hline
case 2&$\cdots ababab \cdots$ &$\in \Q$  &  $\Z^2$& $\Z$\\
\hline
case 3 & $\cdots aabb \cdots$ & $\notin \Q$ & $FG_2$ & $\Z$\\
\hline\hline
\end{tabular}
\end{center}
\end{table}

\subsection{Model sets}
A model set is a point-set obtained by the cut-and-projection scheme. 
There are several formulations of this method and we use here
that of \cite{Moody97}.

Let $H$ be a locally compact abelian group and
$\Lambda\subseteq \RR^d\times H$ a subgroup such that the quotient 
$\RR^d\times H/\Lambda$ is compact. 
In most applications $H=\RR^n$ and $\Lambda$ is a regular lattice. 
Let $\pi:\RR^d\times H\to \RR^d$ be the projection along $H$ and 
$\pi':\RR^d\times H\to H$ the projection along $\RR^d$.  
Let $K\subseteq H$ be a non-empty bounded subset. 
Then
$$\D_K =\{\pi(x)|x\in\Lambda,\pi'(x)\in K\}$$
is a model set with {\em acceptance domain $K$}. 
Note that $\D_K$ depends only on $K'=K\cap\pi'(\Lambda)$. 
We require the following additional assumptions:
 
\begin{description}
\item[{\rm (i)}]
The restriction of $\pi$ to $\Lambda$ is injective.
\item[{\rm (ii)}] The
restriction of $\pi'$ to $\Lambda$ is injective.
\item[{\rm (iii)}] $G = \pi'(\Lambda)$ is dense in $H$.
\item[{\rm (iv)}]  $K$ is the closure of its interior. 
\end{description}
These assumptions are the ones most commonly used and, 
apart from the last one --- which ensures that $\D_K$ is a Delone set ---
they impose no substantial restrictions 
if one wants to describe aperiodic model sets ((i) and (ii) imply that
$\D_K$ has no translational symmetry).

We can construct two semigroups from the data of a cut-and-projection scheme:
the point-set semigroup $\Gamma(\D_K)$, 
which we write more briefly as $\Gamma_K$, 
and the semigroup $\Gamma_K^* =\Gamma(K,\pi'(\Lambda),H)$.
We shall now describe the relationship 
between these two semigroups. 

From conditions (i) and (ii) above,
it follows that $\pi(\Lambda)$ is isomorphic to $\pi'(\Lambda)$. 
We define an isomorphism 
${}^*:\pi(\Lambda)\to \pi'(\Lambda)$ by
$$\pi(x)^* = \pi'(x).$$ 
This allows us to write
$$\D_K = \{y\in \pi(\Lambda) | y^*\in K\}.$$ 
Let us denote by $\C_K$ (resp.\ $\C_K^*$) 
the inverse categories which enter into
the definition of $\Gamma_K$ (resp.\ $\Gamma_K^*$).  
The relevant action on $\C_K$ is that of $\pi(\Lambda)$
by left translation. 
We look at this as a $\Lambda$-action: 
for $g\in\Lambda$ define
$$g\cdot p= \pi(g)+p, \quad p\in \pi(\Lambda).$$
The relevant action on $\C_K^*$ is that of $\pi'(\Lambda)$
by left translation. We look at this as a (right) $\Lambda$-action:
define 
$$g\cdot x = x-\pi'(g), \quad x\in \pi'(\Lambda).$$ 
Note the different sign to the above.

\begin{proposition}
The map 
$\phi:\C_K \to \C_K^*$,
$$\phi(p_2,P,p_1) = (-p_2^*,\bigcap_{x\in
  P}(K-x^*),-p_1^*)$$
is a surjective $\Lambda$-equivariant functor of inverse categories.
\end{proposition}
\bew\ We show first that $\phi$ is well-defined.
We have that a finite $P$ satisfies
$P\subseteq \D_K$ 
iff 
$\forall x\in P$ 
we have that $x^{\ast} \in K$ 
iff 
$0\in\bigcap_{x\in P}(K-x^*)$. 
Define $P^*:=\bigcap_{x\in P}(K-x^*)$.
Thus
$$\phi(p_2,P,p_1) = (-p_{2}^{\ast},P^{\ast},-p_{1}^{\ast}).$$
We now check that
$(-p_{2}^{\ast},P^{\ast},-p_{1}^{\ast}) \in \C_{K}^{\ast}$.
By definition, $p_{i} \in P$
and so $ P^* \subseteq -p_i^*+K$. 
From the definition, $P^{\ast}$ belongs to the underlying semilattice.
Thus $\phi$ is a well-defined function.

Surjectivity is clear from the construction of the model set: any
finite intersection of a shifted acceptance domain which contains $0$
is the image of a pattern $P$ under $*$. 

The identities in $\C_{K}$
are the elements of the form $(p,\{p\},p)$ where $p \in \D_{K}$.
The image of such an element under $\phi$ is
$(-p^{\ast}, K - p^{\ast}  ,-p^{\ast})$ which  
is an identity in $\C_{K}^{\ast}$.
It is easy to check that if $P,Q \subseteq \D_{K}$ then
$(P \cup Q)^{\ast} = P^{\ast} \cap Q^{\ast}$.
Thus
$$\phi(p_2,P\cup Q,q_1) = (-p_2^*,(P\cup Q)^*,-q_1^*) =
(-p_2^*,P^*\cap Q^*,-q_1^*) = (p_2^*,P^*,-p_1^*)(-p_1^*, Q^*,-q_1^*).$$
This shows that $\phi$ is a functor. 

To prove that $\phi$ is $\Lambda$-equivariant,
observe that $\exists g\cdot (p_2,P,p_1)$ iff $P+\pi(g)\subseteq \D$, 
and $\exists g\cdot (q_2,Q,q_1)$ iff $0\in Q-\pi'(g)$. Now
$P+\pi(g)\subseteq \D$ iff $\pi'(g)\in P^*$ is immediate from the definitions.
Finally, if $P+\pi(g)\in\D$ then
$$\phi (g\cdot (p_2,P,p_1)) 
=
\phi (p_2+\pi(g),P+\pi(g),p_1+\pi(g))
=
(-p_2^*-\pi'(g),(P+\pi(g))^*,-p_1^*-\pi'(g))
=
g\cdot(-p_2^*,P^*,-p_1^*),$$
because $(P+\pi(g))^*=P^*-\pi'(g)$. 
This proves $\Lambda$-equivariance.\eb

The above proposition combined with Proposition~2.2.2
implies the following.

\begin{corollary}
$\phi$ induces a surjective homomorphism of semigroups
$[\phi]:\Gamma_K\to \Gamma_K^*$. 
\end{corollary}
\eb

We can obtain an explicit description of the kernel
of the homomorphism $[\phi]$.
The following definition will be useful.
We say two subsets $P,Q\in \D_K$ 
have the same {\em empire} if, for $g\in\Lambda$, 
$\exists g \cdot P \Leftrightarrow \exists g \cdot Q$. 

\begin{lemma} The kernel of $[\phi]$ is the empire congruence $\E$.
\end{lemma}
\bew\ We prove first that 
$\phi(p_2,P,p_1)=\phi(q_2,Q,q_1)$ 
if and only if  
$(p_2,P,p_1) \equiv (q_2,Q,q_1)$,
using the notation introduced prior to Proposition~2.1.3;
this is equivalent to
$p_i=q_i$ and
$P$ and $Q$ have the same empire.

Suppose first that $\phi(p_2,P,p_1)=\phi(q_2,Q,q_1)$. 
It is immediate that $p_i=q_i$ and that $P^{\ast} = Q^{\ast}$.
Now 
$$\exists g+P 
\mbox{ if and only }  
0\in (g+P)^*=P^*-g^*
\mbox{ if and only if }
g^{\ast} \in P^{\ast}.$$
Hence $P^*=Q^*$ implies that  
$\exists g+P$ if and only if $\exists g+Q$.

To prove the converse we make use of our assumptions that $\pi'(\Lambda)$
is dense in $H$, and that $K$ the closure of its interior.
Since $K$ is the closure of its interior,
both $P^*$ and $Q^*$ are closures of their interiors.
Assume $P^*\not\subseteq Q^*$. 
Then $\mbox{int}(P^*) \cap Q^* \neq \mbox{int}(P^*)$ so that
$\mbox{int}(P^*) \backslash \mbox{int}(P^*) \cap Q^*$ is an open subset of
$B =P^*\backslash P^*\cap Q^*$.
Therefore, by the denseness of $\pi'(\Lambda)$,
$\pi'(\Lambda)\cap B$ contains an element $\pi'(g)$. 
It follows that $0\in P^*-\pi'(g)$ but
$0 \notin Q^*-\pi'(g)$. 
Then $\exists \pi(g)+P$ but not $\exists\pi(g) +Q$. 

To finish off the proof, suppose that
$[\phi][p_{2},P,p_{1}] = [\phi][q_{2},Q,q_{1}]$.
Then this is equivalent to
$[\phi (p_{2},P,p_{1})] = [\phi (q_{2},Q,q_{1})]$.
By definition, this implies there exist group
elements $g$ and $h$ such that
$g \cdot \phi (p_{2},P,p_{1})
=
h \cdot \phi (q_{2},Q,q_{1})$.
However, we proved that $\phi$ is $\Lambda$-equivariant.
Thus
$\phi (g \cdot (p_{2},P,p_{1}))
=
\phi (h \cdot (q_{2},Q,q_{1}))$.
It follows that
$[\phi][p_{2},P,p_{1}] = [\phi][q_{2},Q,q_{1}]$
if and only if
there exist 
$(p_{2}',P',p_{1}') \in [p_{2},P,p_{1}]$
and
$(q_{2}',Q',q_{1}') \in [q_{2},Q,q_{1}]$
such that
$\phi (p_{2}',P',p_{1}') 
=
\phi (q_{2}',Q',q_{1}')$.
But by our result above this means precisely that
$(p_{2}',P',p_{1}') \equiv (q_{2}',Q',q_{1}')$.

We have therefore proved that
$[\phi][p_{2},P,p_{1}] = [\phi][q_{2},Q,q_{1}]$
if and only if
$[p_{2},P,p_{1}] \equiv [q_{2},Q,q_{1}]$.\eb

Combining Proposition~3.4.1, Corollary~3.4.2 and Lemma~3.4.3,
we have proved the following result.

\begin{theorem}\label{thm3.4.4} For model sets satisfying 
conditions (i)--(iv) above
$\Gamma_{\D_K}/\emp\cong \Gamma(K,\pi'(\Lambda),H)$.
\end{theorem}
\eb

By Proposition~2.1.3, the empire congruence is idempotent pure;
by Proposition~2.2.1, idempotent pure congruences preserve
universal groups.
Thus $\Gamma_{\D_K}$ has the same universal group as
$\Gamma(K,\pi'(\Lambda),H)$.
If the conditions of Theorem~2.4.7 hold with respect to
the triple $(K,\pi'(\Lambda),H)$
then the universal group of
$\Gamma(K,\pi'(\Lambda),H)$
is isomorphic to $\pi'(\Lambda)$.
We have therefore proved the following theorem,
which is the main result of this paper.

\begin{theorem} 
Consider a model set constructed from the data $(K,G= \pi'(\Lambda),H)$ 
satisfying conditions (i)--(iv) above
and satisfying in addition the following conditions:
 $H$ is connected and simply connected
and $K$ contains an open path-connected subset $V$
such that
$1\in V\subseteq K\subseteq \bigcup_{s\in V} sV$.
Then the universal group of
$\Gamma(\D_K)$ is isomorphic to $G$.
\end{theorem}
\eb

The conditions of the corollary are satisfied for quite a large class
of model sets. In fact, the choice $H=\RR^d$ occurs in many
applications and if the interior of $K$ is connected the above
conditions are met. On the other hand, 
from the example after Theorem~2.4.7 it is not difficult to
construct model sets whose universal group is not isomorphic to
$H_\D$.

\section{Tilings}

Before we discuss in detail the universal groups of connected
semigroups of one-dimensional tilings 
we make some remarks concerning the relation between tilings and
discrete point-sets. We mentioned that from the point of view of
most mathematical theories on aperiodic structures it does not matter
whether one works with tilings or with point-sets, 
the reason being that there are ways of constructing 
point-sets from tilings and vice-versa.

There are various ways of deriving a discrete point-set from a tiling. 
A common choice which is used for polyhedral tilings, i.e.\
tilings whose tiles are polyhedra, 
is to consider the set of vertices as a point-set.
But this is not the only natural choice and, in addition,
it is not invertible, 
because the vertex set does not determine the tiling.  

There are various ways of deriving a tiling from a discrete point-set. 
Two common constructions are the
Voronoi- and the Delone complex associated with the point-set.
These complexes define tilings by declaring the closures of the highest
dimensional cells to be the tiles. The Delone tiling so obtained has
the property that its vertex set corresponds to the discrete point-set
one started with. But not all tilings are Delone tilings of some point-set.

What is important about such constructions is that they are local in
the following sense \cite{Baake}:
given e.g.\ the tiling from which we want to
derive the point-set, such a derivation is called local if the positions
of the points in the point-set which fall
in some ball of radius $r$ in the ambient space $\RR^d$ 
are already determined by 
the tiles of the tiling which lie in a ball (with equal
centre) of possibly larger radius $r'$. 
This idea can be developped into an
equivalence relation, that
of being mutually locally derivable, and a Delone set and its
associated Delone tiling are an example of mutually locally derivable
structures. Furthermore, given a tiling, one can easily improve the 
way described above to derive a Delone set from a tiling to get a set 
which would yield a mutually locally derivable with the tiling. 
The point is that given two mutually locally derivable structures
(tilings or Delone sets)
the topological groupoids constructed from their
inverse semigroup as in \cite{Lawson98} are equivalent
in the sense of \cite{Muhly}.

Mathematical theories which are derived from these groupoids up to
equivalence are
therefore not sensitive to whether one works with tilings
or point-sets. 

However there is one important difference between patterns of tilings
and finite subsets of discrete point-sets: whereas we can make a
distinction on whether a set covered by a pattern is connected or not 
and also tell when a pattern covers a given $r$-ball
this does not work for finite subsets of a
Delone set. 
Therefore, if we
inspect a finite subset of a point set 
we cannot derive a pattern from it 
unless we know that it exhausts all points of the set in some finite
ball. This implies 
that the inverse semigroup of the point-set cannot be expected to be
topologically equivalent \cite{Kellendonk97b} to that of a tiling which
is mutually locally derivable with the set.

\subsection{Tiling semigroups}\setcounter{theorem}{0}

The construction of the tiling semigroup has been discussed in detail
in \cite{Lawson98,KL2000} to which we refer the reader. 

\begin{proposition}
The universal group of the tiling semigroup $\Gamma(\T)$ is
a homomorphic image of that of the connected tiling semigroup $S(\T)$.
\end{proposition}
\bew\ We regard $S(\T)$ as an inverse subsemigroup of $\Gamma (\T)$.
Let 
$\tau:\Gamma(\T)\to G(\Gamma(\T))$
and
$\tau':S(\T)\to G(\Gamma(\T))$
be the respective universal morphisms to the universal groups.
The composition of the inclusion morphism 
$\iota \colon \: S(\T)\hookrightarrow 
\Gamma(\T)$ with the universal morphism $\tau$
yields a morphism  $S(\T)\to G(\Gamma(\T))$,
so that by the universal property there is a homomorphism
$\psi:G(S(\T))\to G(\Gamma(\T))$
such that $\tau \iota = \psi \tau'$.
It follows that $\psi$ maps
$\tau'([p_2,P,p_1])$ to $\tau([p_2,P,p_1])$.

Let $[p_2,P,p_1] \in \Gamma (\T)$.
Then there exists a connected pattern $Q$ 
such that $P \subseteq Q$.
It follows that $[p_{2},Q,p_{1}] \in S(\T)$ 
and $[p_{2},Q,p_{1}] \leq [p_{2},P,p_{1}]$.
Thus $\tau ([p_{2},P,p_{1}]) = \tau ([p_{2},Q,p_{1}])$,
and so $\psi$ maps $\tau' ([p_{2},Q,p_{1}])$
to $\tau ([p_{2},P,p_{1}])$.
But by the remark prior to Proposition~2.2.1,
the image of a universal morphism to a universal group
generates the group.
Thus we have proved that $\psi$ is surjective.\eb

\subsection{Universal groups for 1-dimensional tilings}\setcounter{theorem}{0}

In this section, we shall calculate the universal groups
of semigroups for 1-dimensional tiling,
and show that the connected tiling semigroups are always free.

To handle 1-dimensional tiling semigroups efficiently, 
we need some notation \cite{Kellendonk97a}.
Let $\Sigma$ be an alphabet,
and $\Sigma^{\ast}$ the free monoid on $\Sigma$.
If $u \in \Sigma^{\ast}$ and $u = xyz$ where $x,y,z \in \Sigma^{\ast}$,
then $y$ is said to be a {\em factor} of $u$,
$x$ is a {\em prefix} of $u$ and $z$ is a {\em suffix} of $u$.

A {\em one-dimensional tiling} $\mathcal{T}$ 
is just a bi-infinite string over $\Sigma$:
a function from $\mathbb{Z}$ to $\Sigma$.

With every tiling $\mathcal{T}$ we can associate a $\Sigma$-language
$L(\mathcal{T})$, 
called the {\em language of $\mathcal{T}$},
which consists of all finite, non-empty strings which occur in $\mathcal{T}$.
The language $L(\mathcal{T})$ has the additional property that
if $x \in L(\mathcal{T})$ then all {\em non-empty}
factors of $x$ also belong to $L(\mathcal{T})$.
More generally, 
we say that a language $L$ is {\em factorial}
if it possesses this property.

The connected tiling semigroup $S(\mathcal{T})$ of a tiling $\mathcal{T}$
is actually entirely determined by the language $L(\mathcal{T})$ 
of $\mathcal{T}$. For this reason, we shall extend our range in
this section slightly, by describing the inverse semigroup
$S(L)$ of any factorial language $L$.
In the case where $L$ is the language of a tiling
then $S(L)$ is just the connected tiling semigroup.

Let $L$ be a factorial language over an alphabet $\Sigma$.
The set $S(L)$ is defined as follows.
It consists of all strings over the alphabet
$$\Sigma 
\cup \{\acute{a}\colon \: a \in \Sigma \} 
\cup \{\grave{a} \colon \: a \in \Sigma \} 
\cup \{\check{a} \colon \: a \in \Sigma \}$$
which have the following forms:
\begin{itemize} 

\item $x\check{a}y$ where $x,y \in \Sigma^{\ast}$, $a \in \Sigma$
and $xay \in L$;

\item $u\grave{a}v\acute{b}w$ where $u,v,w \in \Sigma^{\ast}$, $a,b \in \Sigma$
and $uavbw \in L$;

\item $u\acute{a}v\grave{b}w$ where $u,v,w \in \Sigma^{\ast}$, $a,b \in \Sigma$
and $uavbw \in L$.

\end{itemize}
Recall that $\acute{a}$ is the {\em acute accent}
on the letter $a$ and 
and $\grave{a}$ is the {\em grave accent}
on the letter $a$.
The {\em check accent} is to be regarded as simultaneously
grave and acute.
If $p \in S(L)$, then $\delta (p)$ denotes the underlying
string.

If $p \in S(L)$ is a string then the acute accent
in $p$ marks the {\em in-letter}
and the grave accent in $p$ marks the {\em out-letter};
the respective accents are called the {\em in-accent}
and {\em out-accent}.

A product (denoted by $\otimes$)
is defined on $S(L)$ as follows.
Let $p,q \in S(L)$. 
Place $p$ above $q$ so that the 
in-letter of $p$ is above the out-letter of $q$.
We say that $p$ and $q$ {\em match}
if, ignoring accents, 
they agree on their overlap. 
If $p$ and $q$ do not match then $p \otimes q$ is not defined.
If $p$ and $q$ do match then glue the strings together
on their overlap, erasing the in-accent of $p$
and the out-accent of $q$ and carry forward the remaining
two accents;
if as a result a grave and acute accent occur together on a 
letter in the resulting string then rewrite them as a check.
If the resulting string belongs to $S(L)$ then define it
to be $p \otimes q$;
if the resulting string does not belong to $S(L)$
then $p \otimes q$ is undefined.
The proof that $S(L)$ is a well-defined inverse semigroup
is straightforward and left to the reader.

\begin{proposition} Let $L$ be a factorial language.
Then $S(L)$ is strongly $F^{\ast}$-inverse.
\end{proposition}
{\bf Proof }It is evident that $S(L)$ is strongly $E^{\ast}$-unitary,
although officially not covered by our result in
\cite{KL2000}, it is proved in the same way.
Let $$M(L) =
\{\grave{a_{1}} \ldots \acute{a_{m}} \colon \:
a_{1} \ldots a_{m} \in L \} \cup
\{\acute{a_{1}} \ldots \grave{a_{m}} \colon \:
a_{1} \ldots a_{m} \in L \}.$$
It is clear that $M(L)$ is precisely the set of maximal elements of $S(L)$.
An arbitrary element of $S(L)$ has one of three forms:
$x\check{a}y$, $u \grave{a} v \acute{b} w$, 
or $u \acute{a} v \grave{b} w$.
Observe that
$$x \check{a} y \leq \check{a}, \,
u \grave{a} v \acute{b} w \leq \grave{a} v \acute{b}, \,
\mbox{ and }
u \acute{a} v \grave{b} w \leq \acute{a} v \grave{b},$$
where 
$\check{a}, \grave{a} v \acute{b}, \acute{a} v \grave{b} \in M$.
It is evident that these maximal elements are uniquely
determined in each case.
Thus $S(L)$ is $F^{\ast}$-inverse.\qed

Define $C(L) \subseteq S(L)$ by
$$C(L) = \{\grave{a_{1}} \ldots \acute{a_{m}} \colon \:
a_{1} \ldots a_{m} \in L \}.$$
Equip $M(L)$ with the product $\circ$ described prior to Proposition~2.2.3.

\begin{proposition} The universal group of $M(L)$
is isomorphic to the universal group of $C(L)$.
In particular, the universal group of $S(L)$ is isomorphic
to the universal group of $C(L)$.
\end{proposition}
\bew\ We describe the forms taken by the product in $M(L)$
in terms of $C(L)$.
Let $x,y \in C(L)$. 
There are four possible types of product in $M(L)$:
\begin{description}
\item[{\rm (1)}] $x \circ y \in C(L)$.
\item[{\rm (2)}] $x^{-1} \circ y^{-1} = (y \circ x)^{-1} \in M(L)$.
\item[{\rm (3)}] $x \circ y^{-1}$: if 
$\delta (y^{-1})$ is a suffix of $\delta (x)$
then $u \circ y = x$ for some $u \in C(L)$
giving $x \circ y^{-1} = u$ (remembering the definition of $\circ$);
else $\delta (x)$ is a suffix of $\delta (y^{-1})$
and $u \circ x = y$ for some $u \in C(L)$
giving $x \circ y^{-1} = u^{-1}$.
\item[{\rm (4)}] $x^{-1} \circ y$: if $\delta (y)$ is a prefix of
$\delta (x^{-1})$ then $y \circ v = x$ for some $v \in C(L)$
giving $x^{-1} \circ y = v^{-1}$;
else $\delta (x^{-1})$ is a prefix of $\delta (y)$ and
$x \circ v = y$ for some $v \in C(L)$
giving 
$x^{-1} \circ y = v$.
\end{description}
Let $\alpha \colon \: C(L) \rightarrow G$
be any morphism to a group $G$. 
Define
$\alpha^{\ast} \colon \: M(L) \rightarrow G$
by
$$\alpha^{\ast} (x) = \left\{ \begin{array}{ll}
\alpha (x) & \mbox{ if $x \in C(L)$}\\
\alpha (x^{-1})^{-1} & \mbox{ if $x^{-1} \in C(L)$}
\end{array}
\right.
$$
Then our description of the form taken by the product in $M(L)$ now makes
it easy to prove that $\alpha^{\ast}$ is the unique morphism
extending $\alpha$ to $M(L)$.
It is obvious that any morphism from $M(L)$ to a group
restricts to a morphism from $C(L)$ to a group.

If we put these two results together, it is straightforward to check that
the universal group of $C(L)$ is the same as the universal group of $M(L)$.

The proof of the final assertion now follows 
from Proposition~2.2.3.\qed

We can now easily describe the universal group of any $S(L)$.

\begin{theorem} Let $L$ be a factorial language. 
Let $L_{2}$ be the set of all strings of length 2 in $C(L)$.
Then the universal group of $S(L)$ is the free group on $L_{2}$.
\end{theorem}
{\bf Proof }Let $L$ be a factorial language
over the alphabet $\Sigma = \{a_{1}, \ldots, a_{n}\}$.
For convenience, 
we assume that every letter in $\Sigma$ actually occurs in $L$.
By Proposition~4.2.2, 
we need only calculate the universal group of $C(L)$.
Apart from the idempotents $\check{a_{i}}$ for $i = 1,\ldots,n$,
every other element $C(L)^{\ast}$ is a product of elements of length 2, 
because if $w = \grave{b_{1}} \ldots \acute{b_{m}} \in L$
then we can write
$$w = 
(\grave{b_{1}}\acute{b_{2}}) \circ 
(\grave{b_{2}}\acute{b_{3}}) \circ
\ldots 
\circ (\grave{b_{m-1}}\acute{b_{m}}).$$
Since $L$ is factorial, 
each string $\grave{b_{i}}\acute{b_{i+1}} \in L$
for $i = 1,\ldots, m-1$.
Thus $C(L)$ is generated by its idempotents and elements of length 2.

Let $L_{2}$ be the set of all strings of length 2 in $C(L)$,
put $r = |\,L_{2}\,|$,
and let $X_{j}$ be $r$ distinct symbols. 
Let 
$$\phi' \colon \: L_{2} \rightarrow \{X_{j} \colon \: 1 \leq j \leq r\}$$
be a fixed bijection.
Let $FG_{r}$ be the free group on the $r$ symbols $X_{j}$.

Define a function 
$\phi \colon \: C(L)  \rightarrow FG_{r}$ 
as follows.
Let $w \in C(L)$.
If $w$ is an idempotent then define $\phi (w) = 1$;
if $w$ is a string of length 2 in $C(L)$
then define $\phi (w) = \phi' (w)$;
if $w$ is a string of length 3 or more then we can write
it uniquely as a product
$w = w_{1} \circ  \ldots  \circ w_{p}$
of elements $w_{k}$ of length 2, 
in which case define
$\phi (w) = \phi' (w_{1}) \ldots \phi' (w_{p})$.

We show first that $\phi$ is a morphism from $C(L)$ to $FG_{r}$.
Let $u,v \in C(L)$ be arbitrary non-idempotent elements
such that $u \circ v$ is defined.
Then the last letter of $u$ equals the first letter of $v$
and the string underlying $u \circ v$ belongs to $L$. 
Let $u = u_{1} \ldots  u_{s}$ 
and $v = v_{1} \ldots  v_{t}$
where the $u_{i}$ and the $v_{j}$ are strings of length 2. 
Let $u_{s} = \grave{a}\acute{b}$ 
and $v_{1} = \grave{b}\acute{c}$ 
where $a,b,c \in \Sigma$.
It follows that
$u \circ v =  
u_{1} \circ \ldots \circ u_{s} 
\circ 
v_{1} \circ \ldots \circ  v_{t}$
is the correct representation of $u \circ v$ as a
product of strings of length 2.
It is evident that in this case
$\phi (u \circ v) = \phi (u)\phi (v)$.
If now either $u$ or $v$ is an idempotent
such that $u \circ  v$ is defined
then it is easy to see from the definition of $\phi$ 
that $\phi (u \circ v) = \phi (u) \phi (v)$.
It follows that $\phi$ is a morphism.
 
Finally, we show that $\phi$ is the universal morphism to a group.
Let $\theta \colon \: C(L) \rightarrow G$ be any morphism to a group.
For each $u \in L_{2}$, 
let $\theta (u) = g_{j}$
where $\phi'(u) = X_{j}$.
Define $\alpha \colon \: FG_{r} \rightarrow G$
by $\alpha (X_{j}) = g_{j}$ and then extend to the whole
of $FG_{r}$ by freeness.
The functions $\alpha \phi$ and $\theta$
agree on the elements of $L_{2}$ of $C(L)$
and both map non-zero idempotents to the identity of $G$.
Thus $\alpha \phi = \theta$,
and it is clear that $\alpha$
is the unique homomorphism with this property.
Thus $\phi$ is the universal morphism.\qed 

\begin{corollary}
Let $\T$ be a one-dimensional tiling and $n$ be the number of
equivalence classes of consecutive pairs of tiles. The universal group of its
connected semigroup is the free group generated by $n$ elements.
\end{corollary} 

\end{document}